\DeclareSymbolFont{AMSb}{U}{msb}{m}{n}
\documentclass[reqno,centertags,10pt,a4paper,noamsfonts]{amsart}
\pdfoutput=1
\usepackage[svgnames]{xcolor}
\usepackage[english]{babel}
\usepackage[babel=true,expansion=alltext,protrusion=alltext-nott,final]{microtype}
\usepackage[margin={3cm},marginparwidth={2.5cm}]{geometry}
\usepackage[sb]{libertine}
\usepackage[libertine,bigdelims,vvarbb]{newtxmath}
\usepackage[cal=cm]{mathalpha} 

\usepackage{mathtools}
\usepackage[inline]{enumitem}
\numberwithin{equation}{section}

\providecommand{\mr}[1]{\href{http://www.ams.org/mathscinet-getitem?mr=#1}{MR~#1}}
\providecommand{\zbl}[1]{\href{https://zbmath.org/?q=an:#1}{Zbl~#1}}

\providecommand{\doi}[2]{\href{https://doi.org/#1}{#2}}

\newcommand{\C}{\mathcal{C}}

\DeclareMathOperator{\weaklystar}{\rightharpoonup\kern-2.2ex ^* \, \,}
\DeclareMathOperator{\sgn}{sgn}

\newcommand{\R}{\mathbb R}
\newcommand{\N}{\mathbb N}

\renewcommand{\C}{\mathbb C}

\newcommand\norm[1]{\lVert #1 \rVert}

\newcommand\inner[1]{\langle #1 \rangle}

\newcommand\scpr{\boldsymbol{\cdot}}

\newcommand{\ra}{\rightarrow}

\newcommand{\mL}{\mathrm{L}}

\renewcommand{\phi}{\varphi}

\newcommand{\mH}{\mathrm{H}}
\newcommand{\mW}{\mathrm{W}}

\newcommand{\ee}{\mathrm{e}}

\theoremstyle{plain}
\newtheorem{theorem}{Theorem}[section]
\newtheorem{proposition}[theorem]{Proposition}

\newtheorem{lemma}[theorem]{Lemma}

\newtheorem*{theorem*}{Theorem}

\theoremstyle{definition}
\newtheorem{definition}[theorem]{Definition}
\newtheorem{remark}[theorem]{Remark}
\newtheorem*{remark*}{Remark}

\usepackage[pagebackref=false]{hyperref}
\hypersetup{
	linktoc=all,
	bookmarksdepth=3,
	colorlinks=true,
	linkcolor=Green,citecolor=Orange,urlcolor=DarkBlue,
	pdftitle={A non-degeneracy theorem for interacting fermions in one-dimension},
	pdfauthor={Thiago Carvalho Corso}, 
	pdfsubject={},
	pdfkeywords={}} 
\begin{document}
\numberwithin{table}{section}
\title[Non-degeneracy theorem for interacting fermions in one-dimension]{A non-degeneracy theorem for interacting fermions in one dimension}

\author[T.~Carvalho~Corso]{Thiago Carvalho Corso}
\address[T.~Carvalho Corso]{Institute of Applied Analysis and Numerical Simulation, University of Stuttgart, Pfaffenwaldring 57, 70569 Stuttgart, Germany}
\email{thiago.carvalho-corso@mathematik.uni-stuttgart.de}


\thanks{\emph{Funding information}:  DFG -- Project-ID 442047500 -- SFB 1481.  \\[1ex]
\textcopyright 2025 by the authors. Faithful reproduction of this article, in its entirety, by any means is permitted for noncommercial purposes.}
\begin{abstract}
In this paper, we show that the ground-state of many-body Schr\"odinger operators for electrons in one dimension is non-degenerate. More precisely, we consider Schr\"odinger operators of the form
\begin{align*}
    H_N(v,w) = -\Delta + \sum_{i\neq j}^N w(x_i,x_j) + \sum_{j=1}^N v(x_i) \quad \mbox{acting on $\bigwedge^N \mL^2([0,1])$,}
\end{align*}
where the external and interaction potentials $v$ and $w$ belong to a large class of distributions. In this setting, we show that the ground-state of the system with Fermi statistics and local boundary conditions is non-degenerate and does not vanish on a set of positive measure. In the case of periodic and anti-periodic (or more general non-local) boundary conditions, we show that the same result holds whenever the number of particles is odd and even, respectively. This non-degeneracy result seems to be new even for regular potentials $v$ and $w$. As an immediate application of this result, we prove eigenvalue inequalities and the strong unique continuation property for eigenfunctions of the single-particle one-dimensional operators $h(v) = -\Delta +v$. In addition, we prove strict inequalities between the lowest eigenvalues of different self-adjoint realizations of $H_N(v,w)$. 
\end{abstract}
\setcounter{tocdepth}{1}
\maketitle
\tableofcontents
\setcounter{secnumdepth}{2}

\section{Introduction}

The primary goal of this paper is to study the spectral properties of Schr\"odinger operators of the form
\begin{align*}
    H_N(v,w) = -\Delta + \sum_{i \neq j}^N w(x_i,x_j) + \sum_{j=1}^N v(x_j) 
\end{align*}
acting on the space of antisymmetric wave-functions $\wedge^N \mL^2(I)$, where $I = (0,1)$ (or $\R$), and $v$ and $w$ belongs to suitable class of distributions. Such operators are used to describe interacting (spinless) electrons living in one-dimensional space and are therefore ubiquitous in quantum mechanics. 

In the mathematical literature, Schr\"odinger operators with distributional potentials have been extensively studied. Early contributions include \cite{RS72,AGHK84,Bra85,AGHH88,Her89}, while more recent works can be found in \cite{AK00,EGNT15,LR15,EP16,BFK+17}. In particular, special attention has been devoted to the case of one-dimensional Sturm-Liouville operators \cite{Kur97,SS99,AHM08,Dav13,EGNT13a,EGNT13b,MK14,Gul19}. However, most (if not all) of these works focus on the case of what could be called "single-particle electronic systems", i.e., where $H = - \Delta + V$ acts on $\mL^2(\Omega)$ for some $\Omega \subset \R^d$ instead of the antisymmetric $\mL^2$-space. Of course, the spectral properties of the many-particle non-interacting Hamiltonian $H_N(v,0)$ acting on $\wedge^N L^2(I)$ are determined be the spectral properties of the single-particle operator $h(v) = -\Delta +v$ acting on $\mL^2(I)$. Unfortunately, this is no longer true for interacting systems, which is arguably the most interesting case in applications.

In fact, our main motivation for studying the interacting operator $H_N(v,w)$ stems from recent advances \cite{SPR+24,Cor25a} towards establishing a mathematically rigorous formulation of Kohn-Sham Density Functional Theory (DFT) \cite{HK64,KS65}, a highly successful and widely used approach to study the electronic structure of materials. More precisely, the class of distributional potentials studied here appeared in recent developments \cite{SPR+24} concerning the $v$-representability problem, which, roughly speaking, consists in characterizing the set of all possible single-particle ground state densities of operators of the form $H_N(v,w)$ for a given class of potentials $v$. While a detailed discussion on the foundation of DFT and the $v$-representability problem goes beyond the scope of the current paper, we emphasize that the results obtained here are used to provide a completely rigorous formulation of Kohn-Sham DFT for electrons in one-dimension in a companion article \cite{Cor25b}. 

In this work, we focus instead on the ground-state properties of $H_N(v,w)$. More precisely, our main contributions can be stated as follows:
\begin{enumerate}[label=(\roman*)]
\item We prove that the ground-state of self-adjoint realizations of $H_N(v,w)$ with local boundary conditions (BCs) is non-degenerate and almost everywhere non-vanishing. In particular, this applies to the Dirichlet and Neumann realizations of $H_N(v,w)$. \item In the case of non-local boundary conditions such as periodic and anti-periodic BCs, we prove that the ground-state is non-degenerate and almost everywhere non-vanshing whenever the number of particles $N$ is odd and even, respectively.
\item As an application of these non-degeneracy results, we obtain eigenvalue inequalities and prove the strong unique continuation property for eigenfunctions of the one-dimensional single-particle operator $h(v) = -\Delta +v$.
\item In addition, we obtain strict inequalities between the lowest eigenvalues (or ground-state energies) of different self-adjoint realizations of $H_N(v,w)$. \end{enumerate}

\section{Main Results}
    In this section, we state our main results precisely. We then outline the key steps in the proofs and how these steps are organized throughout the paper.
\subsection{Notation} We start with some notation. Throughout this paper, we let $I=(0,1)$ be the open unit interval and set $I_N \coloneqq (0,1)^N$ for any $N\in \N$. 

We denote by $\mH^1(I)$ the Sobolev space of functions $f\in\mL^2(I)$ with weak derivative $\partial_x f \in \mL^2(I)$. More generally, for $1\leq p \leq \infty$ and open $\Omega \subset \R^d$, we denote by $\mW^{1,p}(\Omega)$ the usual Sobolev spaces of functions in $\mL^p(\Omega)$ with weak gradient in $\mL^p(\Omega)$, and by $\mW^{-1,q}(\Omega)$, where $1/q+1/p = 1$, the dual space of $\mW^{1,p}(\Omega)$. In addition, we denote by $\mH^{1/2}(\partial \Omega)$ the standard $1/2$-Sobolev (or Besov) space along the boundary $\partial \Omega$. 

For a given closed subspace $B\subset \mH^{1/2}(\partial \Omega)$, we denote by $H^1_B(\Omega)$ the Sobolev spaces of functions with boundary trace on $B$, i.e.,
\begin{align*}
    \mH^1_B(\Omega) \coloneqq \{ f\in \mH^1(\Omega) : \gamma f \in B\},
\end{align*}
where $\gamma : \mH^1(\Omega) \rightarrow \mH^{1/2}(\partial \Omega)$ is the Dirichlet boundary trace operator. The letter $B$ here stands for boundary conditions.

Next, we define $\mathcal{V}$ as the following space of generalized external potentials:
\begin{align*}
\mathcal{V} \coloneqq \mH^{-1}(I;\R) = \{ v \in \mH^{-1}(I) : v(\phi) \in \R \quad\mbox{for any real-valued $\phi\in \mH^1(I)$}\}.
\end{align*}
Similarly, we define $\mathcal{W}$ as the following space of generalized interaction potentials:
\begin{align*}
    \mathcal{W} &\coloneqq \cup_{q>2} \mW^{-1,q}(I_2;\R) \\
    &= \{w \in \mW^{-1,q}(I_2): \mbox{for some } q>2\quad \mbox{and}\quad w(\phi) \in \R \quad \mbox{for any real-valued $\phi \in \mW^{1,p}(I_2)$}\}.
\end{align*} 

For $N\in \N$, we denote by $\mathcal{H}_N$ the usual space of electronic wave-functions, i.e., the antisymmetric $N$-fold tensor product
\begin{align*}
    \mathcal{H}_N = \bigwedge^N \mL^2(I).
\end{align*}
For $v\in \mathcal{V}$, $w\in \mathcal{W}$, and $B\subset \mH^{1/2}(\partial I_N)$ as above, we denote by $H_N^B(v,w)$ the $N$-particle Hamiltonian
\begin{align*}
    H_N^B(v,w) = -\Delta + \sum_{i\neq j}^N w(x_i,x_j) + \sum_{j=1}^Nv(x_i) \quad \mbox{acting on $\mathcal{H}_N$ with boundary conditions in $B$.} 
\end{align*}
More precisely, $H_N^B(v,w)$ is the unique self-adjoint operator associated to the sesquilinear form 
\begin{align}
    \mathfrak{a}_{v,w}(\Psi,\Phi) \coloneqq \int_{I_N} \overline{\nabla \Psi(x_1,...,x_N)} \scpr \nabla \Phi(x_1,...,x_N) \mathrm{d}x_1...\mathrm{d} x_N + v\left(\rho_{\Psi,\Phi}\right) + w(\rho^{(2)}_{\Psi,\Phi}) \label{eq:quadratic form0}
\end{align}
with form domain 
\begin{align*}
    \mathcal{Q}_N(B) \coloneqq \mathcal{H}_N \cap \mH^1_B(I_N),
\end{align*} 
where $\rho_{\Psi,\Phi}$ is the overlapping single-particle density 
\begin{align}
    \rho_{\Psi,\Phi}(x) \coloneqq N \int_{I_{N-1}} \overline{\Psi(x,x_2,...,x_N)} \Phi(x,x_2,...,x_N) \mathrm{d} x_2...\mathrm{d} x_N , \label{eq:density def}
\end{align}
and $\rho^{(2)}_{\Psi,\Phi}$ is the overlapping pair density
\begin{align}
    \rho^{(2)}_{\Psi,\Phi}(x,y) \coloneqq N(N-1) \int_{I_{N-2}}\overline{\Psi(x,y,x_3,...,x_N)} \Phi(x,y,x_3,...,x_N) \mathrm{d}x_3...\mathrm{d}x_N. \label{eq:pair density def}
\end{align}

\subsection{Main results}
    As a first result, we show that the ground-state of self-adjoint realizations of $H_N(v,w)$ with local boundary conditions is non-degenerate. Precisely, let $\Gamma \subset \partial I_N$ and define
    \begin{align*}
        \mH^1(I_N;\Gamma) = \{ \Psi \in \mH^1(I_N) : (\gamma \Psi)\rvert_\Gamma = 0\}.
    \end{align*}
    Then the following holds.

    \begin{theorem}[Non-degeneracy theorem with local BCs] \label{thm:main} Let $v\in \mathcal{V}$, $w\in \mathcal{W}$, $N\in \N$ and $\Gamma \subset \partial \Omega$. 
    Then the operator $H_{N}(v,w;\Gamma)$ defined as the self-adjoint realization of~\eqref{eq:quadratic form0} with form domain 
    \begin{align*}
        \mathcal{Q}_N(\Gamma) \coloneqq\mH^1(I_N;\Gamma) \cap \mathcal{H}_N
    \end{align*}
    has a unique ground-state $\Psi$ and $\Psi \neq 0$ almost everywhere.
    \end{theorem}

\begin{remark}[Unique continuation] The fact that $\Psi$ does not vanish on a set of positive Lebesgue measure is known as the strong unique continuation property or principle (UCP). While the UCP is known to hold for distributional solutions of the Schr\"odinger equation with a large class of multiplicative potentials, see, e.g. \cite{JK85,Sog90,Wol93,KT01,Hor07,Gar18}, the author is not aware of similar results for many-body Schr\"odinger operators with distributional potentials. In fact, many of the difficulties encountered throughout our proofs can be traced back to this lack of unique continuation results. Interestingly though, Theorem~\ref{thm:main} allow us to establish a unique continuation result for the single-particle operator $h(v) = -\Delta +v$; see Theorem~\ref{thm:non-interacting}.
\end{remark}

The next result shows that Theorem~\ref{thm:main} can be extended to the case of periodic and anti-periodic boundary conditions. In this case, however, we need an additional assumption on the number of particles $N\in \N$. To state this result precisely, let us define, for any $\alpha \in \R\setminus \{0\}$ the space
\begin{align}
    \mH^1_{\alpha}(I_N) \coloneqq \left\{ \Psi \in \mH^1(I_N) :  (\gamma\Psi)\rvert_{I_{k-1}\times\{0\} \times I_{N-k}} - \alpha (\gamma \Psi)\rvert_{I_{k-1} \times \{1\} \times I_{N-k}}= 0\quad \mbox{for any $1\leq k \leq N$} \right\}. \label{eq:non-local BCs}
\end{align}
We say that these spaces have non-local boundary conditions as the boundary values along any face of $\partial I_N$ depend on the values along the opposite face\footnote{We exclude the case $\alpha =0$ for two reasons. First, this case corresponds to local boundary conditions and are therefore covered by Theorem~\ref{thm:main}. Second, the notation $\mH^1_0(I_N)$ is reserved to the case of Dirichlet boundary conditions along the whole boundary.}. Note that $\alpha = 1$ corresponds to periodic BCs and $\alpha = -1$ to antiperiodic BCs. For such spaces, the following result holds.

\begin{theorem}[Non-degeneracy theorem with non-local BCs] \label{thm:non-local} Let $\alpha \in \R\setminus \{0\}$ and suppose that 
\begin{align}
    \alpha (-1)^{N-1} > 0. \label{eq:parity condition}
\end{align}
Then the self-adjoint realization $H_N^\alpha(v,w)$ of~\eqref{eq:quadratic form0} with form domain 
\begin{align*}
    \mathcal{Q}_N^\alpha \coloneqq \mH^1_\alpha(I_N)\cap \mathcal{H}_N
\end{align*}
has a unique ground state $\Psi$ and $\Psi \neq 0$ almost everywhere.
\end{theorem}

\begin{remark}[On the condition on the number of particles] \label{rem:parity condition} Theorem~\ref{thm:non-local} shows that periodic systems are non-degenerate for an odd number of particles, while anti-periodic systems are non-degenerate for an even number of particles. This result is optimal in the sense that there are explicit examples where the ground-state is non-degenerate if and only if condition~\eqref{eq:parity condition} holds, see Appendix~\ref{app:Laplace}. On the other hand, one can construct examples of periodic/anti-periodic systems where the ground-state is non-degenerate for any given number of particles $N$. 
\end{remark}

\begin{remark}[Further extensions] 
   Theorem~\ref{thm:main} can also be extended in the following directions:
    \begin{enumerate}[label=(\arabic*)]
    \item (Real line) Theorem~\ref{thm:main} can be extended to the whole real line $I=\R$ (or half-line). In this case, the existence of a ground-state is not guaranteed and the statement should be changed to "if a ground-state exists" then all of the conclusions of Theorem~\ref{thm:main} holds. 
    \item (Many-body interactions) We can allow a larger class of distributional interaction potentials $w$. In fact, Theorem~\ref{thm:main} holds for $M$-body interactions $w_M$ with arbitrary $M\in \N$, provided that $w_M$ lies in suitable dual Sobolev spaces (see Section~\ref{sec:examples}). Notably, the case of two-body 3D Coulomb interactions, which could be seen as a "six-body" interaction of the form
    \begin{align*}
        w_6(x_1,x_2,x_3,y_1,y_2,y_3) = \left(\sum_{j=1}^3 (x_j-y_j)^2\right)^{-\frac12},
    \end{align*}
    is also admissible. In particular, one can show that the ground-state of the 3D many-particle system with Coulomb interaction, restricted to maximally antisymmetric wave-functions, i.e., wave-functions that are antisymmetric with respect to the exchange of all possible one-dimensional coordinates, is also non-degenerate.
    \item (Elliptic operators) Theorems~\ref{thm:main} and~\ref{thm:non-local} also hold with the Laplacian replaced by more general elliptic operators
    \begin{align*}
        \mathcal{L}_a \Psi \coloneqq \sum_{i,j=1}^N \partial_{x_i} \left(a_{ij}(x) \partial_{x_j} \Psi\right), \quad \mbox{with $a_{ij}(x) \in \mL^\infty(I_N;\R)$.}
    \end{align*}
    The only condition needed for our proof is that weak subsolutions $\mathcal{L}_a u \leq 0$ satisfy the strong maximum principle. 
    \item (More general boundary conditions) We can further extend Theorem~\ref{thm:main} to a larger class of self-adjoint realizations of $H_N(v,w)$. Roughly speaking, the maximal class of spaces for which our proof works is given by all $\mH^1_B(I_N)$ where the restricted space
    \begin{align*}
        \{\Psi \rvert_{S_N} : \Psi \in \mH^1_B(I_N) \cap \mathcal{H}_N \} \quad \mbox{where} \quad S_N = \{(x_1,...,x_N) \in I_N : x_1<x_2... <x_N\}
    \end{align*}
    is invariant under taking the real and positive part, see Theorem~\ref{thm:PFFermi}. In particular, Theorem~\ref{thm:non-local} still holds with the following mixed space of non-local and local boundary conditions
    \begin{align}
        \mH^1_{\alpha}(I_N;\Gamma) = \mH^1_{\alpha}(I_N) \cap \mH^1(I_N;\Gamma), \label{eq:mixed}
    \end{align}
    for $\alpha \in \R\setminus \{0\}$ satisfying~\eqref{eq:parity condition}.
    \end{enumerate}
    \end{remark}
    
Interestingly, Theorems~\ref{thm:main} and~\ref{thm:non-local} can be applied to obtain information on the whole spectrum of single-particle operators. More precisely, let $v\in \mathcal{V}$, $B \subset \C^2$ be a (closed) subspace, and denote by $h_B(v) = -\Delta_B + v$ the self-adjoint operator associated to the form
\begin{align}
    \mathfrak{a}_v(\psi, \psi) = \int_I |\nabla \psi(x)|^2 \mathrm{d} x + v(|\psi|^2) \quad \mbox{with form domain $\mathcal{Q}_1(B) = \{\psi \in \mH^1(I) : (\psi(0), \psi(1)) \in B\}$.} \label{eq:single-particle form} 
\end{align}
Note that $h_B(v) = H_1^B(v,0)$. Then by applying Theorems~\ref{thm:main} and~\ref{thm:non-local} to self-adjoint realizations of the non-interacting operator $H_N(v,0)$, we obtain the following result.
\begin{theorem}[Single-particle operators] \label{thm:non-interacting} Let $v\in \mathcal{V}$ and $h_\alpha(v) \coloneqq h_{B_\alpha}(v)$ where $B_\alpha = \{ \beta =(\beta_0,\beta_1) \in \C^2 : \beta_0 - \alpha \beta_1 = 0\}$ for some $\alpha \in \R\setminus \{0\}$. Then $h_{\alpha}(v)$ has discrete spectrum, all eigenfunctions are almost everywhere non-vanishing, and the eigenvalues $\{\lambda_k\}_{k \in \N}$ ordered in non-decreasing order and counting multiplicity satisfy
\begin{align}
   \begin{dcases} \lambda_{2k-1} < \lambda_{2k}, \quad \mbox{if $\alpha \geq 0$, }\\
   \lambda_{2k} < \lambda_{2k+1}, \quad \mbox{if $\alpha \leq 0$,} \end{dcases}, \quad \mbox{for any $k\in \N$.} \label{eq:single eigenvalue}
\end{align}
Moreover, in the case where $B$ is one of the following sets
\begin{align}
    \{0\} \times \C, \quad \C \times \{0\}, \quad \{(0,0)\}, \quad \mbox{or} \quad \C^2, \label{eq:separable}
\end{align}
the same holds but all the eigenvalues are simple.
\end{theorem}

\begin{remark*}[Sturm-Liouville theory] The result from Theorem~\ref{thm:non-interacting} in the case of local boundary conditions~\eqref{eq:separable} is, at least for regular potentials $v$, well-known in the theory of Sturm-Liouville operators, see, e.g.,  \cite[Theorem 1.9.2]{Zet21}. Nevertheless, the approach presented here seems to be completely different from previous proofs and has the important advantage that it extends to the interacting case. \end{remark*}

Next, we present a monotonicity result for the ground-state energy of $H_N(v,w)$. This result shows that the ground-state energy is strictly monotone with respect to enlarging the Dirichlet set of the form domain of $H_N(v,w)$. Precisely, let us denote by $\Gamma_N$ the symmetrization of the set $\Gamma \subset \partial I_N$, i.e., 
\begin{align*}
    \Gamma_N \coloneqq \cup_{\sigma \in \mathcal{P}_N} \sigma(\Gamma) \quad \mbox{where} \quad \sigma(x_1,...,x_N) = (x_{\sigma(1)}, ..., x_{\sigma(N)} ),
\end{align*}
and $\mathcal{P}_N$ is the set of permutations of $\{1,...,N\}$. Then we have the following result.
\begin{theorem}[Monotonicity with respect to the Dirichlet set] \label{thm:eigenvalues}
Let $\Gamma \subset \Gamma' \subset \partial I_N$ be such that $\Gamma_N' \setminus \Gamma_N$ has non-empty (relative) interior in $\partial I_N$. Let $\lambda_1(\Gamma)$ and $\lambda_1(\Gamma')$ be respectively the lowest eigenvalues of $H_N(v,w;\Gamma)$ and $H_N(v,w;\Gamma')$. Then 
\begin{align*}
    \lambda_1(\Gamma)< \lambda_1(\Gamma').
\end{align*}
Moreover, if $\alpha \in \R\setminus \{0\}$ satisfies~\eqref{eq:parity condition} and $\lambda_1(\Gamma)$ and $\lambda_1(\Gamma')$ denote respectively the lowest eigenvalues of the self-adjoint realizations with form domain $\mathcal{Q}_N^\alpha(\Gamma) \coloneqq \mH^1_\alpha(I_N;\Gamma)\cap \mathcal{H}_N$ respectively $\mathcal{Q}_N^\alpha(\Gamma') =\mH^1_\alpha(I_N;\Gamma')\cap \mathcal{H}_N$, where $\mH^1_\alpha(I_N;\Gamma)$ is defined in~\eqref{eq:mixed}, then the same inequality holds.
\end{theorem}

\begin{remark*}[Strict inequality]\label{rem:strict}
 The important  point in Theorem~\ref{thm:eigenvalues} is the strict inequality, since the weaker inequality $\leq$ trivially follows from the variational principle.
\end{remark*}

\begin{remark*}[Topological assumption]
Our main tool to establish Theorem~\ref{thm:eigenvalues} is a weak unique continuation result along the boundary, namely Theorem~\ref{thm:UCP boundary}. As this result applies only to relatively open subsets of the boundary, the assumption that $\Gamma_N'\setminus \Gamma_N$ has non-empty interior is crucial for our proof. Lifting this assumption seems like a challenging problem. 
\end{remark*}

\subsection{Outline of the proofs}
We now discuss the key steps and main challenges in the proofs of our main results.

The proof of Theorems~\ref{thm:main} and~\ref{thm:non-local} is carried out in Section~\ref{sec:main proof} and relies on two main steps. The first step consists in proving a Perron-Frobenius (PF) theorem for Schr\"odinger operators of the form $-\Delta +V$ for a large class of distributional potentials $V$ (see Theorem~\ref{thm:PF}). In other words, we show that the ground-state of such operators is non-degenerate and strictly positive almost everywhere. This step relies on well-known semigroup techniques and a perturbative argument borrowed from Reed and Simon \cite{RS78}[Section XII.12] (see also \cite{GJ70, Gol77,FS75,Ges84} for similar arguments).  Combining this result with a density result for dual Sobolev spaces (Lemma~\ref{lem:density result}), we can then establish a PF theorem for Schr\"odinger operators of the form $H_N(v,w) = -\Delta + \sum_{j} v_j + \sum_{j\neq i} w_{i,j}$ under a suitable positivity assumption on their form domain. Unfortunately, these results only apply to the operator without Fermi statistics, i.e., without antisymmetry restrictions on the wave-function. Indeed, whenever antisymmetry is present the wave-function must have both positive and negative values and can not be almost everywhere strictly positive. In particular, as highlighted in \cite[page 207]{RS78}, PF results cannot be applied to electronic systems. 

It turns out, however, that this is not entirely true for the case of systems with antisymmetry with respect to exchange of one-dimensional coordinates. Indeed, in the second and key step of our proof, we show that the operator $H_N(v,w)$, acting on antisymmetric functions in $I_N$, is unitarily equivalent to a reduced Schr\"odinger operator acting on functions on the simplex $S_N \coloneqq \{(x_1,...,x_N) \in (0,1)^N : x_1< x_2... < x_N\}$, but without symmetry constraints. This reduction is based on the rather simple but surprisingly useful observation that the box $I_N$ can be tilled by reflections of the simplex $S_N$. This reduction then allow us to apply the aforementioned PF theorem to establish the non-degeneracy and the strong unique continuation property of the ground-state for the electronic Hamiltonian $H_N(v,w)$. At this step, condition~\eqref{eq:parity condition} can be shown to be equivalent to positivity of the form domain of the reduced operator and is therefore crucial for our proof. This step is carried out in Section~\ref{sec:simplex reduction}.

\begin{remark}[Jordan-Wigner transform and Bose-Fermi duality in 1D] \label{rem:Jordan-Wigner} As brought to my attention by an anonymous referee, it turns out that the strategy described above is related to the Jordan-Wigner (JW) transform \cite{JW28} or bosonization approach, which is apparently well-known in the physics community. More precisely, while this transform is usually applied to establish an equivalence between spin systems and 1D lattice fermions (see \cite[Chapter 6]{Gia04}), a continuum analog has been used in a celebrated work of Girardeau \cite{Gir60} to establish the existence of a Fermi-Bose map. This continuum JW transform (or sign map) of Girardeau can be shown to be equivalent to the simplex reduction strategy used here, see Remark~\ref{rem:Bose-Fermi equivalence}. Nevertheless, we emphasize that \cite{Gir60} mostly deal with the case of Bosons with purely hard-core interactions (the Tonks-Girardeau gas), for which the equivalent fermionic system is non-interacting.
\end{remark}

In Section~\ref{sec:non-interacting} we present the proof of Theorem~\ref{thm:non-interacting}. This result is an immediate application of Theorems~\ref{thm:main} and~\ref{thm:non-local}, and a simple lemma that connects the boundary conditions of the non-interacting operator $H_N(v,0)$ to the boundary conditions of the single-particle operator $h(v) = -\Delta +v$.

The proof of Theorem~\ref{thm:eigenvalues} is inspired by ideas from previous works on similar eigenvalue inequalities \cite{Fil05,GM09,BRS18,Roh21}. More precisely, the core idea of the proof is to apply the beautiful domain extension argument of Filonov \cite{Fil05} to establish a unique continuation result along the boundary (cf. Theorem~\ref{thm:UCP boundary}). Here, however, we face two fundamental problems when trying to apply this argument. First, it relies on unique continuation results which are not available in the case of distributional potentials. Second, and perhaps more critical, is the fact that the non-local boundary conditions are such that the faces of $I_N$ are "glued" to each other, and therefore, no extension argument is possible. 

The first problem can be overcome by using the unique continuation property of the ground-state guaranteed by the PF theorem. Note that this property is only proven for the ground-state, which explains why we can only obtain inequalities for the lowest eigenvalue of $H_N(v,w)$. To deal with the second difficulty, we rely on a distributional definition of the Neumann trace (or normal derivative). While this definition is the same as in previous works, see, e.g., \cite{GM09} where the Laplacian with non-local (Robin) boundary conditions is considered, we emphasize that the situation here is more delicate for two reasons. First, the non-local boundary conditions here are imposed directly in the form domain of the operator, as opposed to adding a non-local form along the boundary as done in \cite{GM09}. Second, an eigenfunction $\Psi$ of $H_N(v,w)$ has, in general, no more regularity than $\mH^1$. In particular, its Laplacian has no more regularity than $\mH^{-1}$, and therefore, the usual definition of the weak Neumann trace, which requires $\Delta \Psi$ to be at least in $\mH^{s}$ for $s > -1/2$, see \cite[(3.17)]{GM09}, is not meaningful. 

Fortunately, despite this lack of regularity, we can still show that a meaningful sense can be given to the Neumann trace of eigenfunctions along any open subset of the boundary where said eigenfunction vanish. This result crucially relies on the locality of the operator $H_N(v,w)$, a property that indirectly underlies many of the arguments used here. We then obtain an explicitly formula for this weak Neumann trace that could be new and of independent interest; see Lemma~\ref{lem:Ntrace formula}. By combining this formula with a positivity property of the boundary conditions, which is inherited by the ground-state (see Theorem~\ref{thm:PFFermi}), we can show that the Neumann trace of the ground-state vanishes whenever its Dirichlet trace vanishes. This allows us to extend the form domain of $H_N(v,w)$ and carry out the extension argument of Filonov to obtain a contradiction, thereby completing the proof of the boundary unique continuation result in Theorem~\ref{thm:UCP boundary}. Theorem~\ref{thm:eigenvalues} then follows from standard arguments.

\section{Background on Sobolev Spaces}
\label{sec:background} 

In this section we recall some well-known results about Sobolev spaces that will be used throughout our proofs. 

\subsection{Definitions and basic properties} Let us start by recalling the precise definition of Sobolev spaces and setting-up some notation. 

Throughout this section $\Omega \subset \R^d$ will be a bounded, open, and connected subset of $\R^d$ with Lipschitz boundary. This means that, for any $x\in \partial \Omega$, there exists a neighborhood $x\in U \subset \R^d$ such that, up to a rigid motion (i.e., translation and rotation), $U \cap \partial \Omega$ is the graph of a Lipschitz function. On such domains, the Sobolev spaces can be defined as follows.

\begin{definition}[Sobolev spaces]\label{def:Sobolev} For $1\leq p \leq \infty$, we denote by $\mW^{1,p}(\Omega)$, the space of (complex-valued) functions $f\in \mL^p(\Omega)$ with weak gradient $\nabla f \in \mL^p(\Omega;\C^n)$ endowed with the norm
\begin{align}
    \norm{f}_{1,p}^p \coloneqq\norm{f}_{\mL^p}^p + \norm{\nabla f}_{\mL^p}^p.
\end{align}
Moreover, we denote by $\mW^{1,p}_0(\Omega)$ the closure of the space $C^\infty_c(\Omega)$ with respect to the $\mW^{1,p}$-norm. 
\end{definition}

\begin{definition}[Dual Sobolev spaces] \label{def:dualSobolev} For $1\leq p \leq \infty$, we denote by $\mW^{-1,p}(\Omega)$, the dual space of $\mW^{1,q}(\Omega)$, where $q$ is the H\"older conjugate of $p$ (i.e., $1/p+1/q = 1$). More precisely, 
\begin{align*}
    \mW^{-1,p}(\Omega) \coloneqq \{ T: \mW^{1,q}(\Omega) \rightarrow \C \quad \mbox{linear and continuous}\}
\end{align*}
endowed with the norm
\begin{align*}
    \norm{T}_{-1,p} = \sup_{f\in \mW^{1,q} \setminus\{0\}} \frac{|T(f)|}{\norm{f}_{1,q}}.
\end{align*}
Similarly, we denote by $\mW^{-1,p}_0(\Omega)$ the dual spaces of $\mW^{1,q}_0(\Omega)$. 
\end{definition}

\begin{remark}[Notation]  For $p=2$, we use the standard notation $H^1(\Omega)$ instead of $\mW^{1,2}(\Omega)$. Moreover, whenever clear from the context, we shall omit the domain $\Omega$ and simply write $\mH^1$, $\mW^{1,p}$, $\mW^{1,p}_0$ and so on. 
\end{remark}

For later reference, let us recall the celebrated Gagliardo-Nirenberg-Sobolev (GNS) inequality. We refer, e.g., to \cite[Theorem 12.83]{Leo17} for a proof of the general $1\leq p,q \leq \infty$ case, and to \cite{BM18} for precise conditions on the validity of GNS for Sobolev spaces of fractional order.

\begin{lemma}[GNS interpolation inequality]\label{lem:GNS} Let $\Omega \subset \R^d$ be a bounded open and connected domain with Lipschitz boundary. Then for any $f\in \mW^{1,p}(\Omega)$ with $2\leq p <\infty$ such that $\theta = d/2 - d/p \in [0,1]$ we have 
\begin{align}
    \norm{f}_{\mL^p} \lesssim \norm{f}_{\mH^1}^{\theta} \norm{f}_{\mL^2}^{1-\theta} \quad\mbox{where}\quad \theta = \frac{d}{2} - \frac{d}{p}.\label{eq:GNS}
\end{align}
For $d=1$, the case $p =\infty$ is also allowed.
\end{lemma}

Another fundamental property of Sobolev functions on Lipschitz domains is the existence of a trace along the boundary. Let us recall this fact here as well. For a proof, see \cite[Section 18.4]{Leo17}.

\begin{theorem}[Trace of Sobolev functions]\label{thm:trace} For any open, bounded, and Lipschitz $\Omega \subset \R^d$ with $d\geq 2$, there exists a unique continuous trace operator $\gamma : \mH^1(\Omega) \rightarrow \mH^{\frac12}(\partial \Omega)$ satisfying $\gamma f = f \rvert_{\partial \Omega}$ for $f\in \mH^1(\Omega) \cap C(\overline{\Omega})$. Moreover, the trace operator is surjective and there exists (infinitely many) right inverses, i.e., continuous maps $J:\mH^{\frac12}(\partial \Omega) \rightarrow \mH^1(\Omega)$ such that $\gamma Jf = f$ for any $f\in \mH^{\frac12}(\partial \Omega)$. 
\end{theorem}

\begin{remark}[Fractional Sobolev space on the boundary] Although the precise definition of $\mH^{\frac12}(\partial \Omega)$ will not be relevant to us here, let us state it at least once. We say that a function $f:\partial \Omega \rightarrow \C$ belongs to $\mH^{1/2}(\partial \Omega)$ if $f$ is measurable with respect to the $\mathscr{H}^{d-1}$ Hausdorff measure and satisfies 
\begin{align*}
    \norm{f}_{\mH^{1/2}(\partial \Omega)}^2 = \int_{\partial \Omega} |f(x)|^2 \mathscr{H}^{d-1}(\mathrm{d}x) + \int_{\partial \Omega \times \partial \Omega} \frac{|f(x)-f(y)|^2}{|x-y|^d}  \mathscr{H}^{d-1}(\mathrm{d}x) \mathscr{H}^{d-1}(\mathrm{d}y) < \infty.
\end{align*}
Alternatively, $\mH^{1/2}(\partial \Omega)$ can be defined via localization and the Fourier transform, see \cite[Chapter 3]{McL00}.
\end{remark}

\begin{remark}[One dimensional case] For $\Omega = I = (0,1)$, any function in $\mH^1(I)$ is continuous up to the boundary. In this case, we set $\mH^{1/2}(\partial I) = \C^2$ and the trace operator reduces to $(\gamma f) = (f(0),f(1))\in \C^2$.
\end{remark}

Using the trace operator, we can define the following subspaces of $\mH^1(\Omega)$. 
\begin{definition}[Trace-restricted Sobolev spaces]\label{def:traceSobolev} For any closed subspace $B\subset \mH^{\frac12}(\partial \Omega)$, we denote by $\mH^1_B(\Omega)$ the space of Sobolev functions with trace in $B$, i.e.,
\begin{align*}
    \mH^1_B(\Omega) = \{ f\in \mH^1(\Omega) : \gamma f \in B\},
\end{align*}
where $\gamma$ is the trace operator of Theorem~\ref{thm:trace}.
\end{definition}

From the continuity of the trace operator, it is easy to see that $\mH^1_B(\Omega)$ is a closed subspace of $\mH^1(\Omega)$. Moreover, any right inverse $J: \mH^{\frac12}(\partial \Omega) \rightarrow \mH^1(\Omega)$ of the trace operator defines an isomorphism $\mH^1(\Omega) \cong \mH^1_0(\Omega)\oplus \mH^{1/2}(\partial \Omega)$ via the map
\begin{align*}
    G: \mH^1(\Omega) \rightarrow \mH^1_0(\Omega) \oplus \mH^{1/2}(\partial \Omega), \quad G(f) = (f- J\gamma (f), \gamma (f)).
\end{align*}
From this decomposition, it is easy to see that any closed subspace of $\mH^1(\Omega)$ containing $C^\infty_c(\Omega)$ is given by $\mH^1_B(\Omega)$ for a unique closed subspace $B\subset \mH^{\frac12}(\partial \Omega)$. More precisely, the following holds.
\begin{proposition}[Intermediate spaces between $\mH^1_0$ and $\mH^1$] \label{prop:closed subspaces} Let $X\subset \mH^1(\Omega)$ be a closed subspace such that $C^\infty_c(\Omega) \subset X$, then there exists a unique closed subspace $B\subset \mH^{1/2}(\partial \Omega)$ such that $X = \mH^1_B(\Omega)$.
\end{proposition}

\subsection{Representation and density on dual Sobolev spaces} The following lemma gives a simple but useful representation of the dual Sobolev spaces.
\begin{lemma}[Representation of dual Sobolev space] \label{lem:dual rep} Let $1<p \leq \infty$, then for any $T\in \mW^{-1,p}(\Omega)$, there exists $\alpha \in \C$ and $V \in \mL^p(\Omega;\C^n)$ such that
\begin{align}
    T(f) = \overline{\alpha} \int_\Omega f(x)\mathrm{d} x + \int_\Omega \overline{V}(x) \scpr \nabla f(x) \mathrm{d} x. \label{eq:dual rep}
\end{align}
\end{lemma}

\begin{proof} First, we note that the map 
\begin{align*}
    i:\mW^{1,q}(\Omega) \rightarrow \C \times \mL^q(\Omega;\C^n) \quad i(f) = \left(\int_\Omega f\mathrm{d} x, \nabla f\right)
\end{align*}
is a continuous immersion, i.e., $i$ is continuous, injective and has closed range. Indeed, the continuity is immediate from the definition of $\mW^{1,p}$ and the injectivity follows from the simple fact that
\begin{align*}
    \nabla f - \nabla g =  0 \quad \mbox{implies} \quad f-g=\mathrm{constant} \quad \mbox{for any distributions $f,g\in \mathcal{D}'(\Omega)$},
\end{align*}
To see that $i$ has closed range, note that if $\{(\int f_n, \nabla f_n)\}_{n\in \N}$ is a Cauchy sequence in $\C \times \mL^q(\Omega;\C^d)$, then by Poincare's inequality
\begin{align*}
\lim_{n\ra \infty} \norm{f_n - f_m - \frac{1}{|\Omega|}\int (f_n-f_m)}_{\mL^q} = 0.
\end{align*}
Thus $f_n$ is a Cauchy in $\mW^{1,q}(\Omega)$; in particular, the limit of this sequence exists and belongs to $f\in \mW^{1,q}$ by completeness. Hence $i(f) = (\int f, \nabla f) =  \lim_{n\rightarrow \infty} (\int f_n , \nabla f_n)$.

Now let $T\in \mW^{-1,p}(\Omega)$. As the inverse of $i$ is continuous from $\mathrm{ran}(i)$ to $\mW^{1,q}$ by the closed graph theorem, the functional $T\circ i^{-1}$ belongs to the dual of $\mathrm{ran}(i)$. Since $\mathrm{ran}(i)$ is a closed subspace of $\C \times \mL^q(\Omega;\C^d)$, by the Hahn-Banach theorem there exists a continuous extension of $T\circ i^{-1}$ to $\C \times \mL^q(\Omega;\C^d)$, denoted here by $T'$. Thus, by the Riesz representation theorem in $\mL^q$ spaces (recall that $1\leq q <\infty$ since $1< p\leq \infty$), there exists $\beta \in \C$ and $V \in \mL^p(\Omega; \C^d)$ such that 
\begin{align*}
	T'(\beta, W) = \overline{\alpha} \beta + \int_\Omega \overline{V}(x) \scpr W(x) \mathrm{d} x, \quad \mbox{for any $\beta\in \C$ and $W\in \mL^q(\Omega;\R^d)$.}
\end{align*}
In particular, for $(\beta,W) = (\int f, \nabla f)$, equation~\eqref{eq:dual rep} holds and the proof is complete.
\end{proof}

An important consequence of the above representation is the following density result.
\begin{lemma}[Density of smooth functions]\label{lem:density result} Let $T \in \mW^{-1,p}(\Omega)$ with $1< p <\infty$, then the space of functions $ \{g+\alpha : g\in C^\infty_c(\Omega), \alpha\in \C\}\subset C^\infty(\overline{\Omega})$ is dense in $\mW^{-1,p}(\Omega)$. 
\end{lemma}

\begin{proof} Let $(\alpha, V) \in \C\times \mL^p(\Omega; \C^d)$ be given by Lemma~\ref{lem:dual rep}. As $C^\infty_c(\Omega)$ is dense in $\mL^p(\Omega)$ for any $1\leq p <\infty$, there exists a sequence of vector fields $V_n \in C_c^\infty(\Omega;\C^d)$ such that $\norm{V_n - V}_{\mL^p} \rightarrow 0$ as $n \rightarrow \infty$. In particular, the sequence $g_n \coloneqq \alpha - \mathrm{div}\, V_n \in C^\infty_c(\Omega)+ \C$ satisfies
\begin{align*}
	\left|T(f) - \int_\Omega \overline{g_n(x)} f(x) \mathrm{d} x \right| &= \left|\int_\Omega \overline{V(x)} \scpr \nabla f(x) \mathrm{d} x + \int_\Omega \mathrm{div}\, \overline{V_n}(x) f(x) \mathrm{d} x \right| = \left| \int_\Omega \overline{(V-V_n)}(x) \scpr \nabla f(x) \mathrm{d} x\right| \\
	&\leq \norm{V - V_n}_{\mL^p} \norm{f}_{1,q}, \quad \mbox{for any $f\in \mW^{1,q}$.}
\end{align*}
Hence, $g_n \rightarrow T$ in $\mW^{-1,p}$, which completes the proof.
\end{proof}

\begin{remark}[Dual representation on restricted Sobolev spaces] Since $\mW^{1,q}_B(\Omega)$ is a closed subspace of $\mW^{1,q}(\Omega)$, any functional in $\mW^{-1,p}_B(\Omega)$ can be continuously extended to a functional in $\mW^{-1,p}(\Omega)$ by the Hahn-Banach theorem. In particular, the representation and density results in Lemmas~\ref{lem:dual rep} and~\ref{lem:density result} hold for $\mW^{-1,p}_B(\Omega)$ and $\mW^{-1,p}_0(\Omega)$ as well.
\end{remark}

\subsection{Regularity of reduced densities} We now present a regularity result for reduced densities of wave-functions with finite kinetic energy. This result is a simple application of the GNS, H\"older's, and Minkowski integral inequalities.
\begin{lemma}[Regularity of reduced densities] \label{lem:regularity reduced densities} Let $\Psi,\Phi \in \mH^1(\Omega \times \Omega')$ where $\Omega \subset \R^d$ and $\Omega' \subset \R^{d'}$ are bounded Lipschitz domains. Let $\rho_{\Psi,\Phi}$ denote the overlapping density
\begin{align}
\rho_{\Psi,\Phi}(x) \coloneqq \int_{\Omega'} \overline{\Psi(x,y)} \Phi(x,y) \mathrm{d} y. \label{eq:general density}
\end{align}
Then we have
\begin{align}
	\norm{\rho_{\Psi,\Phi}}_{1,p} \lesssim \norm{\Psi}_{\mH^1} \norm{\Phi}_{\mH^1}^{d-\frac{d}{p}} \norm{\Phi}_{\mL^2}^{1-d+\frac{d}{p}} + \norm{\Phi}_{\mH^1} \norm{\Psi}_{\mH^1}^{d-\frac{d}{p}} \norm{\Psi}_{\mL^2}^{1-d+\frac{d}{p}},  \quad \mbox{for any}\quad\begin{cases} 1\leq p \leq 2, \quad \mbox{if $d=1$,} \\
  1\leq p < 2, \quad \mbox{if $d=2$,} \\
  1 \leq p \leq \frac{d}{d-1}, \quad \mbox{if $d\geq 3$,} \label{eq:reduced Sobolev}
 \end{cases}
\end{align}
where the implicit constant depends on $p$, $\Omega$, and $\Omega'$, but is independent of $\Psi$ and $\Phi$.
\end{lemma}

\begin{proof} First, observe that
\begin{align*}
 \nabla_x \rho_{\Psi,\Phi}(x) = \int_{\Omega'} \overline{\nabla_x \Psi(x,y)} \Phi(x,y)  + \overline{\Psi(x,y)} \nabla_x \Phi(x,y) \mathrm{d} y
\end{align*}
for any smooth $\Psi, \Phi$. Then we can apply Minkowski integral inequality, H\"older's inequality, and the GNS inequality~\eqref{eq:GNS} to find that
\begin{align*}
	\norm{\nabla_x \rho_{\Psi,\Phi}}_{\mL^p} &= \left(\int_\Omega \left| \int_{\Omega'} (\overline{\nabla_x \Psi} \Phi + \overline{\Psi} \nabla_x \Phi) \mathrm{d} y \right|^p \mathrm{d}x \right)^{\frac{1}{p}} \leq \int_{\Omega'} \left( \int_\Omega (|\nabla_x \Phi \Psi| + |\Phi \nabla_x \Psi|)^p \mathrm{d} x\right)^{\frac{1}{p}} \mathrm{d} y\\
	&\leq \int_{\Omega'} \left(\int_\Omega |\nabla_x \Psi|^2 \mathrm{d} x \right)^{\frac{1}{2}} \left(\int_\Omega |\Phi|^{\frac{2p}{2-p}} \mathrm{d} x \right)^{\frac{2-p}{2p}} \mathrm{d} y + \int_{\Omega'} \left(\int_\Omega |\nabla_x \Phi|^2 \mathrm{d} x \right)^{\frac{1}{2}} \left(\int_\Omega |\Psi|^{\frac{2p}{2-p}} \mathrm{d} x \right)^{\frac{2-p}{2p}} \mathrm{d} y \\
    &\lesssim \int_{\Omega'} \norm{\nabla_x \Psi}_{\mL^2_x} \norm{\Phi}_{\mH^1_x}^{\theta} \norm{\Phi}_{\mL^2_x}^{1-\theta} + \norm{\nabla_x \Phi}_{\mL^2_x} \norm{\Psi}_{\mH^1_x}^{\theta} \norm{\Psi}_{\mL^2_x}^{1-\theta} \mathrm{d} y \\
    &\lesssim \norm{\Psi}_{\mH^1} \norm{\Phi}_{\mH^1}^{\theta} \norm{\Phi}_{\mL^2}^{1-\theta} + \norm{\Phi}_{\mH^1} \norm{\Psi}_{\mH^1}^{\theta} \norm{\Psi}_{\mL^2}^{1-\theta},
\end{align*} 
where $\theta = d(\frac{1}{2} - \frac{2-p}{2p}) = d(1-\frac{1}{p})$. As this holds for any smooth $\Psi, \Phi$, the result follows by an approximation argument. The restriction $p<2$ for $d=2$ is necessary in~\eqref{eq:reduced Sobolev} because the GNS inequality~\eqref{eq:GNS} is not valid in the end point case $(p,d) = (\infty,2)$. On the other hand, the case $(p,d)=(\infty,1)$ is allowed in~\eqref{eq:GNS} and therefore $p=2$ in~\eqref{eq:reduced Sobolev} is allowed for $d=1$.
\end{proof}

An important feature of inequality~\eqref{eq:reduced Sobolev} is that only the dimension of the set $\Omega$ is relevant for the regularity of the reduced densities. Moreover, in the case where $d'=0$, i.e., we consider only $\Psi,\Phi \in \mH^1(\Omega)$ and $\rho_{\Psi,\Phi} = \overline{\Psi} \Phi$, estimate~\eqref{eq:reduced Sobolev} reduces to the standard (Sobolev) regularity expected by the product of two $\mH^1$ functions in $
\R^d$. In particular, it implies that the set $\mH^1(\Omega)$ for an interval $\Omega = (a,b) 
\subset \R$ is an algebra of functions. This property was crucial in \cite{Cor25a,SPR+24} but will play no special role here.

\section{Non-degeneracy of the ground-state}\label{sec:main proof}

Our goal in this section is to prove Theorem~\ref{thm:main}.

\subsection{Perron-Frobenius theorem} We start with a Perron-Frobenius theorem for the free Laplacian without particle statistics. To properly state this result, we introduce the following definition.
\begin{definition}(Positivity and reality preserving) \label{def:positivity} Let $B\subset \mH^{\frac12}(\partial \Omega)$ be a closed subspace and $\mH^1_B(\Omega)$ the space introduced in Definition~\ref{def:traceSobolev}. Then we say that 
\begin{enumerate}[label=(\roman*)]
\item $\mH^1_B(\Omega)$ is reality preserving if $\mathrm{Re}\, \Psi \in \mH^1_B(\Omega)$ for any $\Psi\in \mH^1_B(\Omega)$. 
\item $\mH^1_B(\Omega)$ is positivity preserving if $\Psi_+ = \max\{\Psi,0\} \in \mH^1_B(\Omega)$ for any real-valued $\Psi\in \mH^1_B(\Omega)$. 
\end{enumerate}
\end{definition}

\begin{lemma}[Perron-Frobenius for the free Laplacian] \label{lem:PF laplace} Let $\Omega \subset \R^d$ be open, connected, and bounded subset with Lipschitz boundary, and let $\mH^1_B(\Omega)$ be positivity and reality preserving. Then the Laplacian $-\Delta_B$, defined as the unique self-adjoint operator associated to the quadratic form
\begin{align*}
    \mathfrak{a}_0:\mH^1_B(\Omega) \times \mH^1_B(\Omega) \rightarrow \C, \quad \mathfrak{a}_0(\Psi,\Phi) \coloneqq \int_\Omega \overline{\nabla \Psi(x)} \scpr \nabla \Phi(x) \mathrm{d} x,
\end{align*}
has a non-degenerate ground-state and the unique (up to a global phase) normalized ground-state wave-function is strictly positive everywhere in $\Omega$.
\end{lemma}

\begin{remark*}
    The above result appears in various forms throughout the literature \cite[Chapter 6]{Eva10}, and the general form stated above is likely known to experts. Nevertheless, due to the lack of a reference for this precise statement, we provide the simple proof below.
\end{remark*}
\begin{proof}
    The proof is based on the variational principle and the strong maximum principle. So first, we note that $-\Delta_B$ has discrete spectrum and therefore a ground-state $\Psi$ exists. Next, we claim that the ground-state can be taken real-valued. To see this, note that $\overline{\Phi} \in \mH^1_B(\Omega)$ for any $\Phi \in \mH^1_B(\Omega)$ by the reality preserving property. In particular
    \begin{align*}
        \mathfrak{a}_0(\overline{\Psi}, \Phi) = \int_\Omega \nabla \Psi(x) \scpr \nabla \Phi(x) \mathrm{d} x = \mathfrak{a}_0(\overline{\Phi},\Psi) = \lambda \inner{\overline{\Phi},\Psi} = \lambda \inner{\overline{\Psi},\Phi}, \quad \mbox{for any $\Phi \in \mH^1_B(\Omega)$,}
    \end{align*}
    where $\lambda\geq 0$ is the ground-state energy. Consequently, $\overline{\Psi}$ is also a ground-state of $-\Delta_B$, and therefore $2\mathrm{Re} \Psi = \Psi + \overline{\Psi}$ is a ground-state as well.
    
    Thus, assuming that the ground-state is real-valued, we can define $\Psi_+ = \max\{0,\Psi\}$ and $\Psi_- = (-\Psi)_+$. Since $\mH^1_B(\Omega)$ is positivity preserving by assumption, both $\Psi_+$ and $\Psi_-$ belong to $\mH^1_B(\Omega)$. Moreover, since they have disjoint support and $\Psi = \Psi_+-\Psi_-$, we have
    \begin{align*}
        \mathfrak{a}_0(\Psi_+,\Psi_+) = \int_\Omega \overline{\nabla \Psi_+(x)} \scpr \nabla \Psi_+(x) \mathrm{d}x = \int_\Omega \overline{\nabla \Psi(x)} \scpr \nabla \Psi_+(x) \mathrm{d} x = \lambda \inner{\Psi,\Psi_+} = \lambda \norm{\Psi_+}_{\mL^2}^2.
    \end{align*}
    Hence, it follows from the variational principle that
    \begin{align*}
        \lambda = \min_{\Psi \in \mH^1_B\setminus \{0\}} \frac{\mathfrak{a}_0(\Psi,\Psi)}{\norm{\Psi}_{\mL^2}^2} = \frac{\mathfrak{a}_0(\Psi_+,\Psi_+)}{\norm{\Psi_+}_{\mL^2}^2}.
    \end{align*}
    As a consequence, $\Psi_+$ satisfies the Euler-Lagrange (eigenfunction) equation
    \begin{align*}
        \mathfrak{a}_0(\Phi,\Psi_+) = \int_\Omega \overline{\nabla \Phi(x)} \scpr \nabla \Psi(x) = \lambda \inner{\Phi,\Psi_+}, \quad \mbox{for any $\Phi\in C_c^\infty(\Omega)$,}
    \end{align*}
    or equivalently, $(-\Delta -\lambda)\Psi_+ = 0$ in the distributional sense. From standard elliptic regularity (cf. \cite[Section 6.3]{Eva10}), it follows that $\Psi_+ \in C^\infty(\Omega)$. As $\Psi_+\geq 0$ and $\lambda \geq 0$ , $\Psi_+$ is a subsolution of Laplace equation, 
    \begin{align*}
        -\Delta \Psi_{+} \geq 0. 
    \end{align*}
    Therefore, we can apply the strong maximum principle \cite[Theorem 3 in Section 6.4.2]{Eva10} to conclude that either $\Psi_+(x) >0$ for any $x\in \Omega$ or $\Psi_+ = 0$. Hence, either $\Psi_+ = 0$ or $\Psi_- = 0$, which implies $\Psi(x)$ has constant sign everywhere in $\Omega$. Thus any real-valued ground-state of $-\Delta_B$ is strictly positive (up to multiplication by a constant). As a consequence, the ground-state is non-degenerate as there cannot be two strictly positive functions that are mutually orthogonal.
\end{proof}

Our next goal is to combine the previous result with the perturbative approach in Reed and Simon \cite[Section XIII.12]{RS78} to prove a Perron-Frobenius theorem for Schr\"odinger operators with generalized potentials. More precisely, we shall use the following lemma, which is a reformulation of results extracted from Reed and Simon \cite[Theorem XIII.43 and Theorem XIII.45]{RS78}. For the proof of this lemma, we refer to their book.
    \begin{lemma}[Non-degeneracy via perturbation] \label{lem:RS} Let $H$ and $H_0$ be self-adjoint operators on $\mL^2(\Omega)$. Suppose the following holds:
    \begin{enumerate}[label=(\roman*)]
        \item (Positive semigroup) \label{it:positivity preserving semigroup} The heat semigroup $\ee^{-tH_0}$ is positivity preserving (i.e., sends non-negative functions to non-negative functions).
        \item (Unperturbed positive ground-state)\label{it:positive ground-state}  $H_0$ has a non-degenerate and  almost everywhere strictly positive ground-state. 
        \item (Approximation property)\label{it:approximation property} There exists $\{V_n\}_{n\in \N} \subset \mL^\infty(\Omega)$ such that $H_0 + V_n$ converges in the strong resolvent sense to $H$. 
    \end{enumerate}
    Then $H$ has a non-degenerate and almost everywhere strictly positive ground-state (provided that a ground-state exists). 
    \end{lemma}
    
    Next, let us recall the Beurling-Deny criterion to establish positivity of the heat semigroup \cite[Theorem XIII.50]{RS78}. Again, we refer to \cite{RS78} for the proof.

\begin{lemma}[Beurling-Deny criterion] \label{lem:Beurling-Deny} Let $H\geq 0$ be a self-adjoint operator on $\mL^2(\Omega)$ and $\mathfrak{a}_H:\mathcal{Q}(H) \rightarrow \R$ be the associated quadratic form. Then the following are equivalent:
\begin{enumerate}[label=(\alph*)]
\item The heat semigroup $e^{-tH}$ is positivty preserving, i.e., sends non-negative functions to non-negative functions.
\item $e^{-tH}$ is reality preserving, the form domain $\mathcal{Q}(H)$ is positivity preserving, and
\begin{align*}
    \mathfrak{a}_H(u_+) + \mathfrak{a}_H(u_-) \leq \mathfrak{a}_H(u),\quad \mbox{for any real-valued $u\in \mathcal{Q}(H)$.}
\end{align*}
\item $\mathcal{Q}(H)$ is absolute value preserving, i.e., $|\Psi| \in \mathcal{Q}(H)$ for any $\Psi \in \mathcal{Q}(H)$, and
\begin{align*}
    \mathfrak{a}_H(|\Psi|) \leq \mathfrak{a}_H(\Psi),\quad \mbox{for any $\Psi \in \mathcal{Q}(H)$.}
\end{align*}
\end{enumerate}
\end{lemma}

\begin{remark}[Reality and posivity preserving vs. absolute value preseving] It is interesting to note that any closed subspace $X\subset \mH^1(\Omega)$ is positivity and reality preserving if and only if it is absolute value preserving. To be precise, the if direction follows by noticing that 
\begin{align*}
    |\alpha \Psi + |\Psi||= \sqrt{(\alpha \mathrm{Re}\, \Psi + |\Psi|)^2 + \alpha^2 (\mathrm{Im}\, \Psi)^2} = \sqrt{(1+\alpha^2) |\Psi|^2 + 2\alpha \mathrm{Re}\,\Psi |\Psi|} \in X
\end{align*}
and taking the limit $\lim_{\alpha \downarrow 0} \frac{|\alpha \Psi + |\Psi|| - |\Psi|}{\alpha} = \mathrm{Re} \Psi$ to conclude that $X$ is reality preserving. The other direction follows by defining the Laplacian with form domain $X$ and using $(b)\Rightarrow (c)$ in Lemma~\ref{lem:Beurling-Deny}.
\end{remark}
    
We can now combine the previous lemmas to prove the following version of the Perron-Frobenius theorem.
\begin{theorem}[Perron-Frobenius for Schr\"odinger operators with distributional potentials]\label{thm:PF} Let $\Omega \subset \R^m$ and $\mH^1_B(\Omega)$ be as in Lemma~\ref{lem:PF laplace}. Let $V:\mH^1_B(\Omega) \times \mH^1_B(\Omega) \rightarrow \C$ be a sesquilinear form and suppose that there exists a sequence $\{V_n\}_{n \in \N} \subset \mL^\infty(\Omega)$ of real-valued functions such that
\begin{align*}
    \lim_{n \ra \infty} \norm{V-V_n}_{\mH^1 \rightarrow \mH^{-1}} = \lim_{n\rightarrow \infty} \sup_{\Psi,\Phi \in \mH^1_B(\Omega)\setminus \{0\}} \frac{\left|V(\Psi,\Phi) - \int_\Omega V_n(x) \overline{\Psi(x)} \Phi(x) \mathrm{d} x\right|}{\norm{\Psi}_{\mH^1} \norm{\Phi}_{\mH^1}} =0.
\end{align*}
Then there exists a unique self-adjoint operator $H(V) = - \Delta_B +V$ with quadratic form
\begin{align}
    \mathfrak{a}_V(\Psi,\Phi) = \int_\Omega \overline{\nabla \Psi(x)} \scpr \nabla \Phi(x) \mathrm{d} x + V(\Psi, \Phi), \quad \mbox{for $\Phi,\Psi \in \mH^1_B(\Omega)$.} \label{eq:H(V) form}
\end{align}
Moreover, the ground-state of $H(V)$ is non-degenerate and the unique (up to a global phase) ground-state wave-function can be chosen strictly positive almost everywhere in $\Omega$.
\end{theorem}

\begin{proof}
    We first show that $V$ is symmetric and $\Delta$-bounded with relative bound $0$. For this, note that, for any $\epsilon>0$ there exists $n \in \N$ such that $\norm{V - V_n}_{\mH^1\rightarrow \mH^{-1}} \leq \epsilon$. Hence,
    \begin{align*}
        \left|V(\Psi,\Phi) \right| \leq \epsilon \norm{\Psi}_{\mH^1} \norm{\Phi}_{\mH^1} + \norm{V_n}_{\mL^\infty} \norm{\Psi}_{\mL^2} \norm{\Phi}_{\mL^2}, \quad \mbox{for any $\Psi, \Phi \in \mH^1_B(\Omega)$,}
    \end{align*}
    and
    \begin{align*}
        |V(\Psi,\Phi) - \overline{V(\Phi,\Psi)}| &= \left|\left(V(\Psi,\Phi) - \int_{\Omega} V_n(x) \overline{\Psi}(x) \Phi(x) \mathrm{d} x\right) + \left(\overline{V(\Phi,\Psi)} - \int_\Omega \overline{V_n(x) \overline{\Phi(x)} \Psi(x)} \mathrm{d} x\right) \right|\\
        &\leq 2\epsilon \norm{\Psi}_{\mH^1} \norm{\Phi}_{\mH^1}.
    \end{align*}
    As $\epsilon>0$ can be taken arbitrarily small, we conclude that $V$ is symmetric and $\Delta$-bounded with relative bound $0$. Hence, by the KLMN theorem (cf. \cite[Theorem X.17]{RS75}), there exists a unique self-adjoint operator with quadratic form~\eqref{eq:H(V) form}. Moreover, as the quadratic form domain of $H(V)$ is $\mH^1_B(\Omega)$, which is compactly embedded in $\mL^2(\Omega)$, the operator $H(V)$ has compact resolvent and therefore purely discrete spectrum. In particular, there exists at least one ground-state. 
    
    To show that the ground-state is non-degenerate and strictly positive almost everywhere, we can now use Lemma~\ref{lem:RS}. For this, we set $H_0 \coloneqq - \Delta_B$, and note that $H_0$ has a non-degenerate and strictly positive ground-state by Lemma~\ref{lem:PF laplace}. Moreover, since $H_0$ is a real operator, i.e., commutes with complex conjugation, the semigroup $e^{-tH_0}$ is reality preserving. Hence, one can easily verify that criterion (b) in Lemma~\ref{lem:Beurling-Deny} holds for $H_0$, and therefore, $\ee^{-tH_0}$ is also positivity preserving. In particular, assumptions~\ref{it:positivity preserving semigroup} and~\ref{it:positive ground-state} from Lemma~\ref{lem:RS} hold. To verify assumption~\ref{it:approximation property}, note that, since $V_n \rightarrow V$ in $\mathcal{B}(\mH^1;\mH^{-1})$, the operators $H_0 + V_n$ converge to $H(V)$ in the strong (in fact norm) resolvent sense. Indeed, this follows from the formula 
\begin{multline*}
    (z-H(V))^{-1} - (z- H_0-V_n)^{-1} \\
    = (z-H(V))^{-1/2} \left(I - \left(I+(z-H(V))^{-1/2}(V-V_n) (z-H(V))^{-1/2}\right)^{-1}\right)(z-H(V))^{-1/2}.
\end{multline*}
(See \cite[Theorem VIII.25.(c)]{RS80} for the detailed argument.) We can thus apply Lemma~\ref{lem:RS} to complete the proof.
\end{proof}
\subsection{Unitary reduction to the simplex} \label{sec:simplex reduction} Throughout this section, we denote by $S_N$ the open simplex
\begin{align*}
    S_N \coloneqq \{(x_1,...,x_N) \in I_N : 0< x_1 < x_2<... <x_N <1\}.
\end{align*}
The key observation that allow us to apply Theorem~\ref{thm:PF} to the operator $H_N^B(v,w)$ is that the hypercube $I_N$ can be decomposed in a disjoint union of reflections of the simplex $S_N$. To make this precise, let us introduce some additional notation. 
\begin{figure}[ht!]
    \centering
    \includegraphics[scale=0.2]{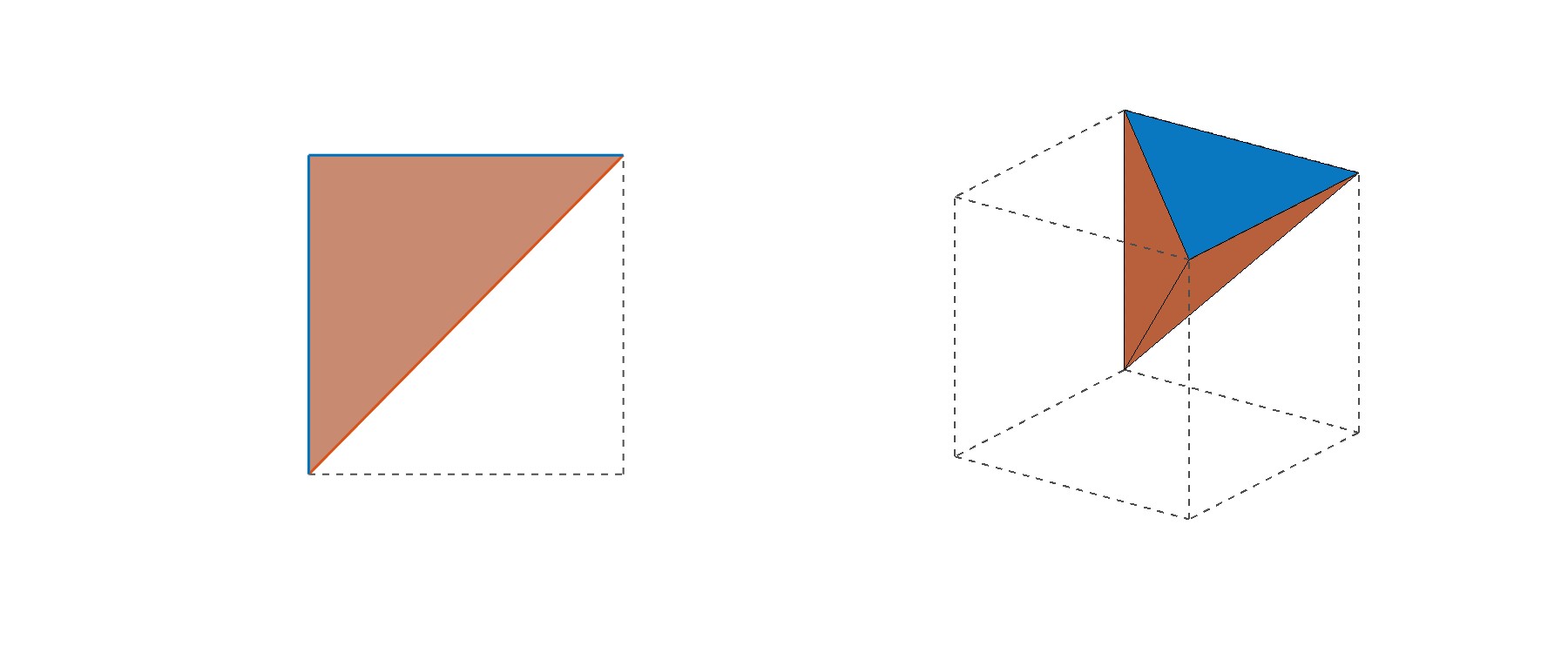}
    \caption{Simplex $S_N$ in the case $N=2$ (left) and $N=3$ (right) with edges of $\partial I_N$ (dashed lines) and interior boundary $\Gamma_{\rm int} \subset \partial S_N$ (in red), where functions in $\mathcal{Q}_N(B;S_N)$ vanish.}
    
 \label{fig:simplex}
\end{figure}

First, we denote by $\mathcal{P}_N$ the set of all permutations of $\{1,...,N\}$, i.e.,
\begin{align*}
    \mathcal{P}_N=\{\sigma: \{1,...,N\} \rightarrow \{1,...,N\}  \mbox{ bijective}\}.
\end{align*}
With some abuse of notation, we also denote by $\sigma$ the linear map of permutation of coordinates
\begin{align}
   \sigma:\R^N \rightarrow \R^N, \quad  \sigma(x_1,...,x_N) = (x_{\sigma(1)},...,x_{\sigma(N)}). \label{eq:permutation}
\end{align}
Recall that the sign of the permutation $\sigma$ is defined as $\mathrm{sgn}(\sigma) = \det \sigma \in \{-1,1\}$. Moreover, we denote by $\{F_j\}_{j=0}^N$ the faces of the boundary of the simplex $S_N$, i..e, 
\begin{align}
    F_j \coloneqq \{ (x_1,...,x_N) \in \overline{S_N} : x_j= x_{j+1} \}, \label{eq:faces def}
\end{align}
with the convention $x_0 =0$ and $x_{N+1} =1$. We can now define the interior part of the boundary of $S_N$ as
\begin{align} 
    \Gamma_{\rm int} \coloneqq \partial S_N \setminus (F_0 \cup F_N) = \{(x_1,..,x_N) \in \overline{S_N} : x_j = x_{j+1} \quad \mbox{for some $j$ and} \quad 0<x_1\leq x_N < 1 \}. \label{eq:interior exterior boundary}
\end{align}
For a visual illustration of $S_N$ and the interior boundary $\Gamma_{\rm int}$, see Figure~\ref{fig:simplex}. 

Let us also define the pushforward map $\sigma_\#$ as
\begin{align*}
    \sigma_\#: \mH^1(S_N) \rightarrow \mH^1(\sigma(S_N)), \quad (\sigma_\#\Psi)(x_1,...,x_N) = \Psi\left(\sigma^{-1}(x_1,...,x_N)\right) = \Psi\left(x_{\sigma^{-1}(1)},...,x_{\sigma^{-1}(N)}\right).
\end{align*}
Finally, let us define the reduced trace-restricted space on $S_N$ as
\begin{align}
    \mathcal{Q}_N(B;S_N) \coloneqq \{\Psi\rvert_{S_N} : \Psi \in \mH^1_B(I_N) \cap \mathcal{H}_N\}. \label{eq:reduced form space}
\end{align} 
Note that any $\Psi \in \mH^1(I_N)\cap \mathcal{H}_N$ has a well-defined and identically zero (Dirichlet) trace along the coincidence hyperplanes $\{x\in I_N: x_j = x_k \mbox{ for $j\neq k$}\}$ by antisymmetry. In particular, any $\Psi\in \mathcal{Q}_N(B;S_N)$ vanishes along the interior boundary $\Gamma^{\rm int}$.

Using this observation, we then have the following key decomposition of $\mathcal{Q}_N(B) =\mathcal{H}_N \cap \mH^1_B(I_N)$.
\begin{lemma}[Unitary reduction to the simplex] \label{lem:simplex reduction} The map $T : \mathcal{Q}_N(B;S_N) \rightarrow  \mathcal{Q}_N(B)$ defined via
\begin{align}
    T(\Psi) = \frac{1}{\sqrt{N!}} \sum_{\sigma \in \mathcal{P}_N} \mathrm{sgn}(\sigma) \sigma_\#(\Psi) \label{eq:extension map}
\end{align}
is an isometric isomorphism.
\end{lemma}

\begin{proof} The key observations here are the following:
\begin{enumerate}[label=(\roman*)]
    \item(Tessellating cover)\label{it:tesselate} The sets $\sigma(\overline{S_N})$ with $\sigma \in \mathcal{P}_N$ tessellate $I_N$, i.e.,
    \begin{align} I_N \subset \cup_{\sigma \in \mathcal{P}_N} \sigma(\overline{S_N})\quad \mbox{and}\quad \sigma(S_N) \cap \tau(S_N) = \emptyset \quad \mbox{ for any $\sigma\neq \tau \in \mathcal{P}_N$.} \label{eq:disjoint property}
\end{align}
 \item(Interior boundary)
 \label{it:boundary} The interior boundary $\Gamma_{\rm int}$ defined in~\eqref{eq:interior exterior boundary} satisfies $\Gamma_{\rm int} = \partial S_N \cap I_N$.
\end{enumerate}
The fact that $\cup_{\sigma}\sigma(\overline{S_N})$ covers $I_N$ is clear.  Indeed, for any $(x_1,...,x_N) \in I_N$ there exists $\sigma \in \mathcal{P}_N$ such that $0\leq x_{\sigma(1)}\leq x_{\sigma(2)} ... \leq x_{\sigma(N)} \leq 1$ and therefore $(x_1,...,x_N) \in \sigma^{-1}(\overline{S_N})$. 

To prove the disjoint property, we first note that, since any $\sigma \in \mathcal{P}_N$ is a homeomorphism in $\R^N$, 
\begin{align*}
    \sigma(S_N) \cap \tau(S_N) = \emptyset \quad \mbox{if and only if} \quad (\tau^{-1}\circ \sigma)(S_N) \cap S_N = \emptyset.
\end{align*}
Hence, it suffices to show~\eqref{eq:disjoint property} for $\sigma \in \mathcal{P}_N$ and $\tau = I$ (identity map). For this, we define
\begin{align*}
    m(\sigma) \coloneqq \max \{j: \sigma^{-1}(k) \leq k \quad \mbox{for any $k\leq j$} \}\cup\{0\}.
\end{align*}
Then, if $\sigma \neq I$ we must have $0\leq m(\sigma)<N$. In particular $\sigma^{-1}(m(\sigma)+1) > m(\sigma)+1 = \sigma^{-1}(m(\sigma)+k)$ for some $k\geq 2$. Then, on the one hand, any $(x_1,...,x_N) \in \sigma(S_N)$ satisfies $x_{\sigma^{-1}(j)} < x_{\sigma^{-1}(j+1)}$ for any $j$ and, in particular, $x_{\sigma^{-1}(m(\sigma)+1)}<x_{\sigma^{-1}(m(\sigma)+k)} = x_{m(\sigma)+1}$. On the other hand, any $(x_1,...,x_N) \in S_N$ satisfies $x_j < x_{j+1}$ for any $j$. Since $\sigma^{-1}(m(\sigma)+1) > m(\sigma)+1$, this implies, in particular, that $ x_{m(\sigma)+1}<x_{\sigma^{-1}(m(\sigma)+1)}$. Therefore, we can not have $(x_1,...,x_N)\in S_N\cap \sigma(S_N)$, which completes the proof of~\ref{it:tesselate}.

We now turn to~\ref{it:boundary}. First, it is clear from the definition~\eqref{eq:faces def} that $F_0\cup F_N \subset \partial I_N$ and therefore $(F_0 \cup F_N )\cap I_N = \emptyset$. Hence, $\partial S_N \cap I_N \subset \Gamma_{\rm int}$. On the other hand, if $(x_1,..,x_N) \in \Gamma_{\rm int} = \partial S_N \setminus (F_0 \cup F_N)$ we must have $0<x_0 \leq x_1 ... \leq x_N < 1$. As every element in the boundary of $I_N$ must have at least one coordinate with value $0$ or $1$, we must have $(x_1,...,x_N) \in I_N$ and therefore item~\ref{it:boundary} holds.

We now show that $T$ is an isometry from $\mathcal{Q}_N(B;S_N)$ to $\mathcal{Q}_N(B)$. To this end, first note that, since $\Psi$ is the restriction of an antisymmetric function, it follows from the definition of $\Gamma_{\rm int}$ (see~\eqref{eq:interior exterior boundary}) that $\gamma \Psi \rvert_{\Gamma_{\rm int}} =0$. Consequently, from property~\ref{it:boundary}, the extension of $\Psi$ to $I_N$ by zero on $I_N\setminus S_N$ belongs to $\mH^1(I_N)$. Hence $T\Psi \in \mH^1(I_N)$. Moreover, if we denote the permutation of $x_j$ and $x_{j+1}$ by $\sigma_j$, then by construction we have
\begin{align*}
    (\sigma_j)_\#(T\Psi) = \frac{1}{\sqrt{N!}} \sum_{\sigma \in \mathcal{P}_N} \sgn(\sigma) (\sigma_j \circ \sigma)_\# (\Psi) = -\frac{1}{\sqrt{N!}} \sum_{\sigma_j \circ \sigma \in \mathcal{P}_N} \mathrm{sgn}(\sigma_j\circ\sigma) (\sigma_j\circ\sigma)_\#(\Psi) = - T\Psi.
\end{align*}
Therefore, $T\Psi \in \mathcal{H}_N \cap \mH^1(I_N)$. Furthermore, since each $\sigma(S_N)$ is disjoint and $\sigma$ an isometry, we have
\begin{align*}
    \norm{T\Psi}_{\mL^2(\Omega)}^2 + \norm{\nabla T\Psi}_{\mL^2(\Omega)}^2 = \frac{1}{N!} \sum_{\sigma \in \mathcal{P}_N} \int_{\sigma(S_N)}  |\Psi(\sigma^{-1} x)|^2 +|\nabla \Psi(\sigma^{-1} x)|^2 \mathrm{d} x = \norm{\Psi}_{\mH^1(S_N)}^2,
\end{align*}
and therefore $T$ is an isometry. Moreover, a similar calculation shows that the restriction map $\Psi \mapsto T^{-1}\Psi = \sqrt{N!} \Psi\rvert_{S_N}$ is an isometry from $\mH^1(I_N) \cap \mathcal{H}_N$ to $\mH^1(S_N;\Gamma_{\rm int})$. As $T^{-1}$ maps $\mathcal{Q}_N(B) = \mH^1_B(I_N) \cap \mathcal{H}_N$ to $\mathcal{Q}_N(B;S_N)$ by definition, it follows that $T:\mathcal{Q}_N(B;S_N) \rightarrow \mathcal{Q}_N(B)$ is surjective, which completes the proof. 
\end{proof}

As the operator $T$ in~\eqref{eq:extension map} is an isometry from $\mathcal{Q}_N(B;S_N)$ to $\mathcal{Q}_N(B)$, the form
\begin{align}
   \widetilde{\mathfrak{a}}_{v,w}(\Psi,\Phi) &\coloneqq   (T^\#\mathfrak{a}_{v,w})(\Psi,\Phi) = \mathfrak{a}_{v,w}(T\Psi, T\Phi) =\int_{I_N} \overline{\nabla (T\Psi)(x)} \scpr \nabla (T\Phi)(x) \mathrm{d} x + v(\rho_{T\Psi,T\Phi}) + w(\rho^{(2)}_{T\Phi,T\Phi}), \nonumber \\
    &= \int_{S_N} \overline{\nabla \Psi(x)} \scpr \nabla \Phi(x) \mathrm{d} x + v(\rho_{T\Psi,T\Phi}) + w(\rho^{(2)}_{T\Phi,T\Phi}), \quad \mbox{for $\Psi,\Phi\in \mathcal{Q}_N(B;S_N)$,}\label{eq:reduced quadratic form}
\end{align}
where the density and pair density are defined as in~\eqref{eq:density def} and~\eqref{eq:pair density def}, is a semibounded closed form on $\mL^2(S_N)$. Hence, we can define the operator 
\begin{align}
    \widetilde{H}^B_N(v,w) \coloneqq T^\#\left(H_N^B(v,w)\right)  = T^{-1} H_N^B(v,w) T\label{eq:reduced H}
\end{align}
as the unique self-adjoint operator associated to the form $\widetilde{\mathfrak{a}}_{v,w}$ with form domain $\mathcal{Q}_N(B;S_N)$. Moreover, this operator is unitarily equivalent to $H^B_N(v,w)$ and we can apply Theorem~\ref{thm:PF} to obtain the following result. 

\begin{theorem}[Perron-Frobenius with Fermi statistics] \label{thm:PFFermi} Let $v\in \mathcal{V}$, $w\in \mathcal{W}$, and $B\subset \mH^{1/2}(\partial I_N)$ be a closed subspace such that the space defined in~\eqref{eq:reduced form space},
\begin{align*}
    \mathcal{Q}_N(B;S_N) = \{ \Psi \rvert_{S_N} : \Psi \in \mH^1_B(I_N) \cap \mathcal{H}_N \},
\end{align*}
is reality and positivity preserving. Then the self-adjoint operator $H_N^B(v,w)$ defined via~\eqref{eq:quadratic form0} (see also Lemma~\ref{lem:operator def}) has a non-degenerate ground-state and the unique (up to a global phase) ground-state wave-function $\Psi$ satisfies $\Psi(x) >0$ for a.e. $x\in S_N$.
\end{theorem}

\begin{proof} By construction, the ground-state of $H^B_N(v,w)$ is given by $T\Psi$ where $\Psi$ is the ground-state of $\widetilde{H}^B_N(v,w)$. Hence, it suffices to show that $\widetilde{H}_N^B(v,w)$ has a non-degenerate and strictly positive ground-state. For this, note that $C^\infty_c(S_N) \subset \mathcal{Q}_N(B;S_N)$. In particular, by Proposition~\ref{prop:closed subspaces}, the form domain $\mathcal{Q}_N(B;S_N)$ is of the form $\mH^1_{B'}(S_N)$ for a suitable closed subspace $B'\subset \mH^{\frac12}(\partial S_N)$. Hence, it suffices to show that $\widetilde{H}^B_N(v,w)$ satisfy the assumptions from Theorem~\ref{thm:PF}. 

To this end, we first use Lemma~\ref{lem:density result} to obtain a sequence of bounded (actually smooth) functions $v_n \in \mL^\infty(I)$ and $w_n \in \mL^\infty(I_2)$ such that $\lim_{n \rightarrow \infty} \norm{v_n -v}_{\mH^{-1}} \rightarrow 0$ and $\lim_{n \rightarrow \infty} \norm{w_n - w}_{\mW^{-1,p}}=0$. Then, by estimate~\eqref{eq:reduced Sobolev} we have
    \begin{align*}
        \lim_{n \rightarrow \infty} \sup_{\Psi,\Phi \in \mH^1(I^N)\setminus\{0\}} \frac{|v_n(\rho_{\Psi,\Phi}) + w_n(\rho^{(2)}_{\Psi,\Phi}) - v(\rho_{\Psi,\Phi}) - w(\rho^{(2)}_{\Psi,\Phi})|}{\norm{\Psi}_{\mH^1} \norm{\Phi}_{\mH^1}} \leq \lim_{n \rightarrow \infty} \norm{v_n-v}_{-1,2} + \norm{w_n-w}_{-1,p} = 0.
    \end{align*}
    As $T: \mathcal{Q}_N(B;S_N) \rightarrow \mathcal{Q}_N(B)$ is continuous, the forms $T^\# v_n$ and $T^\# w_n$, defined as
    \begin{align*}
        (T^\# v_n)(\Psi, \Phi) = v_n(\rho_{T\Psi, T\Phi}) \quad \mbox{and}\quad (T^\# w_n)(\Psi,\Phi) = w_n(\rho^{(2)}_{T\Psi,T\Phi}),
    \end{align*}
    also converge to the forms $T^\#v$ and $T^\#w$ (defined analogously) in the $\mH^1$-form norm, i.e.,
    \begin{align}
        \lim_{n \rightarrow \infty} \sup_{\Psi,\Phi \in \mathcal{Q}_N(B;S_N)\setminus\{0\}} \frac{|(T^\#v_n + T^\# w_n - T^\#v - T^\#w)(\Psi,\Phi)|}{\norm{\Psi}_{\mH^1} \norm{\Phi}_{\mH^1}} = 0. \label{eq:approx property}
    \end{align}
    Moreover, it follows from the definition of $T$ and the tessellating property~\eqref{eq:disjoint property} that, for any $V \in \mL^\infty(I_N)$,
    \begin{align*}
        (T^\# V)(\Psi,\Phi) &= V(T\Psi, T\Phi) = \int_{I_N} V(x) \overline{(T\Psi)(x)} (T\Phi)(x) \mathrm{d} x\\
        &= \frac{1}{N!} \sum_{\sigma,\tau  \in \mathcal{P}_N} \mathrm{sgn}(\sigma) \mathrm{sgn}(\tau) \int_{\sigma(S_N)\cap \tau(S_N)}  \!\!\!\!\!\!\!\!\!\!\!\!\!\!\!\!\!\!V(x) \overline{\Psi(\sigma^{-1} x)} \Phi(\tau^{-1} x) \mathrm{d} x = \frac{1}{N!} \sum_{\sigma \in  \mathcal{P}_N} \int_{S_N} V(\sigma x) \overline{\Psi(x)} \Phi(x) \mathrm{d} x\\
        &= \int_{S_N} \widetilde{V}(x) \overline{\Psi(x)} \Phi(x) \mathrm{d}x , \quad \mbox{where} \quad \widetilde{V}(x_1,...,x_N) = \frac{1}{N!} \sum_{\sigma \in \mathcal{P}_N} V(\sigma x).
    \end{align*}
    In particular, this implies that the pull-back of any form generated by a bounded potential in $I_N$ is a form generated by a bounded potential in $S_N$. As the forms $v_n$ and $w_n$ are generated by the bounded potentials
    \begin{align*}
        V_n(x_1,...,x_N) = \sum_{j=1}^N v_n(x_j) \quad \mbox{and}\quad W_n(x_1,...,x_N) = \sum_{i\neq j} w_n(x_i,x_j), 
    \end{align*}
    the forms $T^\#(v_n + w_n)$ are also generated by bounded potentials. Thus, by estimate~\eqref{eq:approx property}, we have found a sequence of bounded multiplicative potentials approximating $T^\#v + T^\#w$. We can therefore apply Theorem~\ref{thm:PF} to complete the proof.
\end{proof}

\begin{remark}[Bose-Fermi equivalence in 1D] \label{rem:Bose-Fermi equivalence} By considering the symmetric map 
\begin{align*}
    T_+: \mathcal{Q}_N(B;S_N) \rightarrow \mH^1(I_N) \cap \mathcal{H}^{\rm Bos}_N,\qquad \qquad T_{+}(\Psi) \coloneqq \frac{1}{\sqrt{N!}} \sum_{\sigma \in \mathcal{P}_N} \sigma_\#(\Psi),
\end{align*}
where $\mathcal{H}^{\rm Bos}_N = \otimes^N_{\rm sym} \mL^2(I)$ is the symmetric tensor-product, we obtain an isomorphism between $\mathcal{Q}_N(B;S_N)$ and a subset of bosonic (symmetric) wavefunctions vanishing along the hyperplanes $\{x\in I_N: x_j = x_k \quad \mbox{for some $j\neq k$}\}$. This vanishing behaviour can be understood as the effect of a a limiting hard-core pairwise interaction potential $w(x,y) = w_a(x-y) = +\infty$ for $|x-y| \leq a$ and zero otherwise, in the case where $a\downarrow 0$. In particular, we note that the composition map 
\begin{align*}
    T_+ \circ T^{-1}(\Psi)(x_1,...,x_N) = \prod_{j >k} \mathrm{sign}(x_j-x_k) \Psi(x_1,...,x_N)
\end{align*}
is precisely the isomorphism (sign map) between fermions and hard-core bosons explicitly constructed by Girardeau in \cite{Gir60}. Therefore, we can recover the Fermi-Bose equivalence due to \cite{Gir60} for a general class of interactions.
\end{remark}

\subsection{Proofs of Theorems~\ref{thm:main} and~\ref{thm:non-local}} We now prove the following result, which together with Theorem~\ref{thm:PFFermi} completes the proof of Theorems~\ref{thm:main} and~\ref{thm:non-local}. 

\begin{lemma}[Positivity preserving spaces]\label{lem:positivespaces} Let $N \in \N$ and $\Gamma \subset \partial I_N$, then the space
\begin{align}
    \mathcal{Q}_N(\Gamma;S_N) \coloneqq \{\Psi \rvert_{S_N} : \Psi \in \mH^1(I_N;\Gamma)\cap \mathcal{H}_N\} \subset \mH^1(S_N) \label{eq:space1}
\end{align}
is reality and positivity preserving. Moreover, if $\alpha \in \R\setminus \{0\}$ satisfies~\eqref{eq:parity condition}, then the space
\begin{align}
    \mathcal{Q}_N^\alpha(\Gamma;S_N) \coloneqq \{\Psi \rvert_{S_N} : \Psi \in \mH^1_\alpha(I_N;\Gamma) \cap \mathcal{H}_N\} \subset \mH^1(S_N), \label{eq:space2}
\end{align}
where $\mH^1_\alpha(I_N;\Gamma)$ is defined in~\eqref{eq:mixed}, is also reality and positivity preserving.
\end{lemma}

\begin{proof} That both spaces are reality preserving is clear since $\mH^1_\alpha(I_N;\Gamma)$ with $\alpha \in \R$ and $\mathcal{H}_N$ are both reality preserving. For the positivity preserving part, we first deal with $\mathcal{Q}_N(\Gamma;S_N)$. 

In this case, the proof follows from two observations. First, by antisymmetry, the boundary trace of any function $\Psi$ in $\mathcal{Q}_N$ is anti-symmetric and must therefore vanish on the (possibly larger) symmetrized set $\Gamma_N = \cup_{\sigma \in \mathcal{P}_N} \sigma(\Gamma)$. In particular, the restriction $\Psi \rvert_{S_N}$ vanishes on the intersection $\Gamma_N \cap \partial S_N$. The second observation is that $(\gamma \Psi)_+ = \gamma (\Psi_+)$ because this holds for continuous functions and taking the positive part is (e.g., by \cite[Theorem 1]{MM79}) a continuous operation in $\mH^1$. Therefore, $(\Psi\rvert_{S_N})_+ \in \mH^1(S_N;\Gamma_N \cap \partial I_N)$. In particular, the boundary trace of the extension $T(\Psi_+ \rvert_{S_N})$ vanishes on the symmetrized set $\Gamma_N$ and therefore on $\Gamma$. Thus $T(\Psi_+ \rvert_{S_N})$ belongs to $\mH^1(I_N;\Gamma) \cap\mathcal{H}_N$ and satisfies $T(\Psi_+\rvert_{S_N}) \rvert_{S_N} = \frac{1}{\sqrt{N!}} \Psi_+ \rvert_{S_N}$. This shows that $\Psi_+$ belongs to $\mathcal{Q}_N(\Gamma;S_N)$, which implies that this space is positivity preserving. 

For the second statement, it suffices to show that
\begin{align*}
    \mathcal{Q}^\alpha_N(S_N) \coloneqq \{ \Psi \rvert_{S_N} : \Psi\in \mH^1_\alpha(I_N) \cap \mathcal{H}_N \}
\end{align*}
is positivity preserving. Indeed, this follows from the facts that the intersection of two positivity preserving spaces is positivity preserving, and
\begin{align*}
    \mathcal{Q}^\alpha_N(\Gamma;S_N) = \mathcal{Q}^\alpha_N(S_N) \cap \mathcal{Q}_N(\Gamma;S_N).
\end{align*}
To prove that $\mathcal{Q}^\alpha_N(S_N)$ is positivity preserving, we first claim that
\begin{align}
   \mathcal{Q}^\alpha_N(S_N) = \{\Psi \in \mH^1(S_N) :  \Psi\rvert_{\Gamma_{\rm int}} = 0 \mbox{ and } \Psi(0,x') - \alpha (-1)^{N-1}\Psi(x',1) = 0 \quad \mbox{for a.e. $x'\in S_{N-1}$}\}. \label{eq:claim}
\end{align}
Assuming this claim for the moment, the result follows from the following observation. Since $\alpha (-1)^{N-1} \geq 0$ and $\Psi(0,x') =(-1)^{N-1} \alpha \Psi(x',1)$ for (a.e.) $x'\in S_{N-1}$, we have
\begin{align*}
    \Psi(0,x')_+ &=  \max\{ \Psi(0,x'), 0 \} = \max\{ \alpha(-1)^{N-1} \Psi(x',1),0\} = (-1)^{N-1} \alpha \Psi(x',1)_+
\end{align*}
and therefore $\mathcal{Q}^\alpha_N(S_N)$ is positivity preserving. 

To prove the claim, we first note that, by the definition of $\mH^1_{\alpha}(I_N)$ and antisymmetry, any $\Psi \in \mathcal{Q}_N^\alpha$ satisfies
\begin{align*}
    \Psi(0,x_1,...,x_{N-1}) = \alpha \Psi(1,x_1,...,x_{N-1}) = -\alpha \Psi(x_1,1,x_2,...,x_{N-1}) = (-1)^{N-1} \alpha \Psi(x_1,...,x_{N-1},1) 
\end{align*}
for any $x' = (x_1,...,x_{N-1}) \in I_{N-1}$. This proves the inclusion $\subset$ in~\eqref{eq:claim}. To prove the opposite inclusion, we observe that for any $\Psi \in \mH^1(S_N)$ satisfying $\Psi\rvert_{\Gamma_{\rm int}}=0$ and $\Psi(0,x') = (-1)^{N-1} \alpha \Psi(x',1)$ for a.e. $x'\in S_{N-1}$, 
\begin{align}
    (T\Psi) (0,\sigma x') &= \sgn(\sigma) \Psi(0,x') =  \sgn(\sigma) (-1)^{N-1}\alpha \Psi(x',1) = (-1)^{N-1} \alpha (T\Psi)(\sigma x', 1) \nonumber \\
    &= \alpha (T\Psi)(1,\sigma x') \label{eq:bc condition}
\end{align}
for any $\sigma \in \mathcal{P}_{N-1}$ and a.e. $x'\in S_{N-1}$. As $\{\sigma(S_{N-1})\}_{\sigma \in \mathcal{P}_{N-1}}$ covers $I_{N-1}$ up to finitely many hyperplanes of dimension $N-2$, equation~\eqref{eq:bc condition} holds for almost every $x'\in I_{N-1}$. Thus $T\Psi \in \mH^1_\alpha(I_N)$ and satisfies $(T\Psi)\rvert_{S_N} = \frac{1}{\sqrt{N!}} \Psi$, which proves the inclusion $\supset$ in~\eqref{eq:claim} and completes the proof. 
\end{proof}

\subsection{Proof of Theorem~\ref{thm:non-interacting}} \label{sec:non-interacting}

We now prove Theorem~\ref{thm:non-interacting}.

\begin{proof}[Proof of Theorem~\ref{thm:non-interacting}] Let $\alpha \in \R\setminus \{0\}$ and $v\in \mathcal{V}$, then by Lemma~\ref{lem:single-to-many}, the eigenfunctions of the operator $H_N^\alpha(v,0)$ (defined in Theorem~\ref{thm:non-local}) are given by linear combination of Slater determinants of the eigenfunctions of the single-particle operator $h_\alpha(v) \coloneqq H_1^\alpha(v,0)$. In particular, if we denote by $\{\psi_j\}_{j \in \N}$ and $\{\lambda_j\}_{j\in \N}$ the eigenfunctions and eigenvalues of the operator $h_\alpha(v)$ ordered in non-decreasing order, then $H_N(v,0)$ has a non-degenerate ground-state if and only if $\lambda_N < \lambda_{N+1}$. Thus by Theorem~\ref{thm:non-local}, we have $\lambda_N < \lambda_{N+1}$ for any $N\in\N$ satisfying \eqref{eq:parity condition}, in particular, for $N=2k-1$ if $\alpha \geq 0$ and $N=2k$ for $\alpha \leq 0$. This proves~\eqref{eq:single eigenvalue}. The cases of local (or separable) boundary conditions in~\eqref{eq:separable} follows from the same argument. 

To see that the eigenfunctions are almost everywhere non-vanishing, let $E \coloneqq \{\psi_j = 0\}$ and pick $N > j$ satisfying~\eqref{eq:parity condition}. Then the unique ground-state of $H_N^\alpha(v,0)$ is given by
\begin{align*}
    \Psi(x_1,...,x_N) = \frac{1}{\sqrt{N!}} \sum_{\sigma \in \mathcal{P}} \mathrm{sgn}(\sigma) \psi_1(x_{\sigma(1)}) .... \psi_j(x_{\sigma(j)}) ... \psi_N(x_{\sigma(N)}).
\end{align*}
In particular, for $x \in E^N$ we have $x_{\sigma(j)} \in E$ for any $\sigma \in \mathcal{P}_N$, and therefore $\Psi$ vanishes on the set $E^N$. By the strong UCP in Theorems~\ref{thm:main} and \ref{thm:non-local}, this implies that $|E^N| = |E|^N = 0$, which completes the proof.
\end{proof}

\section{Eigenvalue inequalities}

In this section, our goal is to prove Theorem~\ref{thm:eigenvalues}. 

\subsection{Unique continuation along the boundary} \label{sec:UCPboundary} The key ingredient in the proof of Theorem~\ref{thm:eigenvalues} is the following (weak) unique continuation result along the boundary.

\begin{theorem}[Weak unique continuation along the boundary] \label{thm:UCP boundary} Let $\Gamma \subset \partial I_N$, and let $\Psi$ be the unique ground-state of the self-adjoint realization of $\mathfrak{a}_{v,w}$ with form domain $\mathcal{Q}_N(\Gamma)= \mH^1(I_N;\Gamma)\cap \mathcal{H}_N$. Then $\gamma\Psi$ can not vanish identically on a relatively open set 
\begin{align*}
    U\subset \partial I_N\setminus \Gamma_N, \quad\mbox{where}\quad  \Gamma_N = \cup_{\sigma \in \mathcal{P}_N} \sigma(\Gamma).
\end{align*}
Moreover, if $\alpha \in \R\setminus \{0\}$ satisfies~\eqref{eq:parity condition}, then the same holds for the ground-state of the self-adjoint realization with form domain $\mathcal{Q}_N^\alpha(\Gamma) = \mH^1_\alpha(I_N;\Gamma)\cap \mathcal{H}_N$.
\end{theorem}

\begin{figure}[ht!]
    \centering
    \includegraphics[scale=0.2]{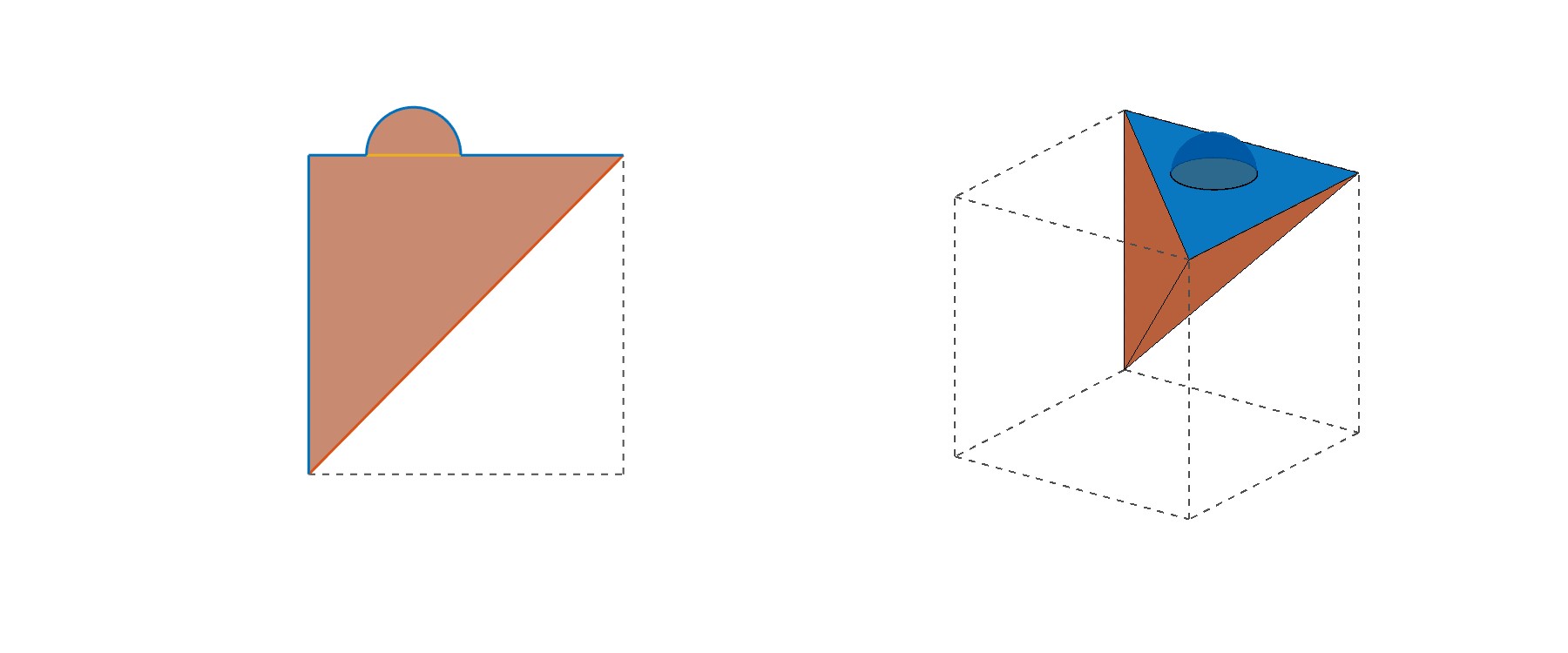}
    \caption{Example of extended set $\Omega$ in the case $N=2$ (left) and $N=3$ (right) with interior boundary $\Gamma_{\rm int}$ (in red), exterior boundary $\partial \Omega \setminus \Gamma_{\rm int}$ (in blue), and $B\cap \partial S_N \subset F_N$ (in yellow).}
    
 \label{fig:simplex_ball}
\end{figure}

\begin{proof}[Proof of Theorem~\ref{thm:UCP boundary}: case of local BCs] To simplify the notation, we set 
\begin{align*}
    \Gamma_S \coloneqq  \Gamma_{\rm int} \cup (\Gamma_N \cap \partial S_N),
\end{align*}
where $\Gamma_{\rm int}$ is the interior boundary of the simplex $S_N$ (see~\eqref{eq:interior exterior boundary}). 
Then, we recall that the ground-state of $H_N(v,w;\Gamma)$ is given by $T\widetilde{\Psi}$, where $T$ is given by~\eqref{eq:extension map} and $\widetilde{\Psi}$ is the ground-state of the reduced operator $\widetilde{H}_N(v,w;\Gamma)$ defined via~\eqref{eq:reduced quadratic form} with form domain~\eqref{eq:space1}. Hence, it suffices to show that $\widetilde{\Psi}$ can not vanish on a relatively open subset of $\partial S_N\setminus \Gamma_S$. To prove this, we shall assume that $\widetilde{\Psi}$ vanishes on a relatively open set $U\subset \partial S_N \setminus \Gamma_N$ and argue by contradiction. 

First, since $U \subset \partial S_N$ is relatively open, we can find an open ball $B\subset\R^N$ such that $B\cap \partial S_N$ is contained in the intersection of $U$ with the interior of one of the exterior faces of $S_N$. So without loss of generality, let us assume that $B \cap \partial S_N \subset U \cap F_0$, where we recall that $F_0 = \{(x_1,...,x_N) \in \partial S_N : x_1 = 0\}$. Next, we define $\Omega\coloneqq S_N \cup B$. Then note that $S_N \subset \Omega \subset \R^N$, $\Gamma_S \subset \partial \Omega$, and $\Omega$ is an open bounded and connected subset of $\R^N$ with Lipschitz boundary (see Figure~\ref{fig:simplex_ball}). In particular, the restriction map 
\begin{align*}
    R: \mH^1(\Omega;\Gamma_S) \rightarrow \mH^1(S_N;\Gamma_S), \quad \Phi \mapsto R \Phi = \Phi \rvert_{S_N},
\end{align*}
is continuous. We can thus define the extended self-adjoint operator $\widetilde{H}_\Omega$ via the sesquilinear form
    \begin{align}
        \widetilde{\mathfrak{a}}_{v,w,\Omega}(\Psi,\Phi) &\coloneqq \int_{\Omega\setminus S_N} \overline{\nabla \Psi(x)} \scpr \nabla \Phi(x) \mathrm{d} x + \widetilde{\mathfrak{a}}_{v,w}(R\Psi,R\Phi) \nonumber \\
        &=\int_{\Omega} \overline{\nabla \Psi(x)} \scpr \nabla \Phi(x) \mathrm{d} x + (T^\#v)(R\Psi,R\Phi) + (T^\# w)(R\Psi, R\Phi), \quad \Psi,\Phi \in \mH^1(\Omega;\Gamma_S).\label{eq:form Omega}
    \end{align}
    Since $(\gamma \widetilde{\Psi}) \rvert_{\Omega \cap \partial s_N} = 0$, the extension by zero of $\widetilde{\Psi}$ to $\Omega$, denoted also by $\widetilde{\Psi}$, belongs to $\mH^1(\Omega;\Gamma_S)$ and satisfies
    \begin{align*}
        \widetilde{\mathfrak{a}}_{v,w,\Omega}(\widetilde{\Psi}, \Phi) = \int_{\Omega} \overline{\nabla \widetilde{\Psi}}\scpr \nabla \Phi \mathrm{d} x + (T^\# (v + w))(R\widetilde{\Psi}, R\Phi) =\widetilde{\mathfrak{a}}_{v,w}(\widetilde{\Psi}, R \Phi) = \lambda \inner{\widetilde{\Psi}, R\Phi} = \lambda \inner{\widetilde{\Psi}, \Phi},
    \end{align*}
    where $\lambda$ is the ground-state energy of $\widetilde{H}_N(v,w;\Gamma)$. In particular, $\widetilde{\Psi}$ is an eigenfunction of $\widetilde{H}_\Omega$. 

    Next, notice that the pullback $R^\# V$ of any bounded multiplicative function $V \in \mL^\infty(S_N)$ is simply the extension of $V$ by zero to $\Omega$, so in particular, a bounded function in $\mL^\infty(\Omega)$. Consequently, the pullback via $R$ of any sequence of bounded functions in $S_N$ approximating the quadratic form $T^\#(v+w)$ in $\mH^1$ (as in~\eqref{eq:approx property}) yields a sequence of bounded functions approximating $R^\# T^\#(v+w)$. Consequently, the operator $\widetilde{H}_\Omega$ defined via~\eqref{eq:form Omega} satisfies the assumptions of Theorem~\ref{thm:PF}, and therefore, has a unique and a.e. strictly positive ground-state $\Psi_\Omega$. 
    
    This now yields a contradiction for the following reason. As we have shown that $\widetilde{\Psi}$ is an eigenfunction of $\widetilde{H}_\Omega$, it must be either orthogonal or proportional to $\Psi_\Omega$. However, $\widetilde{\Psi}$ cannot be orthogonal to $\Psi_\Omega$ as they are both non-negative (up to a global phase), and $\widetilde{\Psi}$ cannot be parallel to $\Psi_\Omega$ as the former vanishes on $\Omega \setminus S_N$ while the later is (a.e.) strictly positive in $\Omega$. This completes the proof.
\end{proof}

In the case of non-local boundary conditions, the proof is more involved. The reason is that the exterior faces $F_0$ and $F_N$ of the simplex are not independent but "glued" to each other in the sense that the boundary conditions along one exterior face must match (up to a constant) the boundary condition along the other exterior face. Hence, we cannot immediately use the extension argument from the previous proof.

To overcome this issue, the first step is to show that, whenever an eigenfunction vanishes on an open subset of the boundary, a meaningful sense can be given to its normal derivative along this set. For this, let us introduce the following definition.

\begin{definition}[Weak Neumann trace] Let $\alpha \in \R\setminus \{0\}$ and $\Psi \in \mathcal{Q}_N^\alpha(\Gamma)$ be an eigenfunction of $H_N^\alpha(v,w;\Gamma)$ with eigenvalue $\lambda$. Then for any $f\in \mH^{1/2}(\partial I_N)$, we define the Neumann trace $\partial_\nu \Psi$ as
\begin{align}
    (\partial_\nu \Psi)(f) = \int_{I_N} \overline{\nabla \Psi(x)} \scpr \nabla F(x) + v(\rho_{\Psi, P_N F}) + w(\rho^{(2)}_{\Psi, P_N F}) - \lambda \inner{\Psi,F} \label{eq:Ntrace}
\end{align}
for any $F\in \mH^1(I_N)$ satisfying $\gamma F = f$, where $\gamma$ is the standard (Dirichlet) trace and $P_N: \mL^2(I_N) \rightarrow \mathcal{H}_N$ is the projection on the antisymmetric space,
\begin{align}
    (P_N F)(x) =  \frac{1}{N!}\sum_{\sigma \in \mathcal{P}_N}\mathrm{sgn}(\sigma) F(\sigma x) \label{eq:antisymmetric projection}
\end{align}
\end{definition}

It is not hard to see that $\partial_\nu \Psi(f)$ is independent of the extension of $f$. Indeed, if we have two $\mH^1$-extensions $F$ and $F'$, then $F-F' \in \mH^1_0(I_N)$, and therefore $P_N (F-F') \in \mH^1_0(I_N) \cap\mathcal{H}_N$; in particular, $P_N(F-F')$ belongs to the form domain $\mathcal{Q}_N^\alpha(\Gamma)$, and consequently,  
\begin{align*}
    &\int_{I_N} \overline{\nabla \Psi(x)} \scpr \nabla(F-F')(x) + v(\rho_{\Psi, P_N(F-F')}) + w(\rho^{(2)}_{\Psi,P_N(F-F')}) - \lambda \inner{\Psi, F-F'} \\
    &=  \int_{I_N} \frac{1}{N!}\sum_{\sigma \in \mathcal{P}_N}\mathrm{sgn}(\sigma)\overline{\nabla (\sigma_\#\Psi)(x)} \scpr \nabla (F-F')(x) \mathrm{d} x + v(\rho_{\Psi, P_N(F-F')}) + w(\rho^{(2)}_{\Psi,P_N(F-F')}) - \lambda \inner{\Psi, P_N(F-F')}\\
    &= \mathfrak{a}_{v,w}(\Psi,P_N(F-F')) - \lambda\inner{\Psi,P_N(F-F')} = 0,
\end{align*}
where the first equality hold because $P_N\Psi = \Psi$ and the second equality because $\Psi$ is an eigenfunction of $H_N^\alpha(v,w;\Gamma)$ with eigenvalue $\lambda$. Moreover, the Neumann trace $f \in \mH^{1/2}(\partial I_N)\rightarrow (\partial_\nu \Psi)(f)$ is continuous since we can pick $F = Jf$ for any right inverse $J : \mH^{1/2}(\partial I_N) \rightarrow \mH^1(I_N)$ of the Dirichlet trace operator $\gamma$. However, it should be noted that this is a completely adhoc definition of the Neumann trace as, a priori, there is no connection between $\partial_\nu \Psi$ and the normal derivative of $\Psi$ along the boundary.

Fortunately, it turns out that this definition is meaningful along the vanishing set of $\Psi$ on the boundary. To make this statement precise, we shall use the following lemma.
\begin{lemma}[Interior approximation] \label{lem:H10 approx} Let $\Psi \in \mH^1(I_N)$ satisfy $(\gamma \Psi)\rvert_U = 0$ for some relatively open subset $U \subset \partial I_N$. Let $V \subset \subset U$ and define
\begin{align*}
    \eta_\epsilon(x) \coloneqq \eta(d_V(x)/\epsilon), \quad \mbox{where}\quad d_V(x) \coloneqq \inf_{y \in V} ||x-y||,
\end{align*}
for some $\eta \in \mathrm{Lip}(\R,\R)$ satisfying $0 \leq \eta \leq 1$, $\eta(x) = 1$ for $|x| \leq 1/2$, and $\eta(x) = 0$ for $|x|\geq 1$. Then the function $\Psi_\epsilon = \Psi (1-\eta_\epsilon) \in \mH^1(I_N;U)$ satisfies
\begin{align*}
    \lim_{\epsilon \downarrow 0} \norm{\Psi - \Psi_\epsilon}_{\mH^1} = \lim_{\epsilon \downarrow 0^+} \norm{\eta_\epsilon \Psi}_{\mH^1} = 0.
\end{align*}
\end{lemma}

\begin{proof}
    Since $\eta_\epsilon(x) \rightarrow 0$ as $\epsilon \downarrow 0$ for every $x\in I_N$, by dominated convergence we have 
    \begin{align*}
        \norm{\Psi-\Psi_\epsilon}_{\mL^2}^2 = \int_{I_N} |\Psi(x)|^2 |\eta_\epsilon(x)|^2 \mathrm{d} x  \rightarrow 0 \quad \mbox{as}\quad \epsilon \downarrow 0.
    \end{align*}
    Next, notice that the map $x\mapsto d_V(x)$ is globally Lipschitz with Lipschitz constant smaller than $1$. Thus from Rademacher's theorem, the function $d_V(x)$ has a bounded weak gradient, and therefore,
    \begin{align*}
        \norm{\nabla(\Psi-\Psi_\epsilon)}_{\mL^2}^2 &\leq 2 \int_{I_N} |\nabla \Psi(x)|^2 |\eta_\epsilon(x)|^2 + |\Psi(x)|^2 \frac{|\dot{\eta}(d_V(x)/\epsilon) |\nabla d_V(x)|^2}{\epsilon^2} \mathrm{d} x \\
        &\lesssim \int_{I_N} |\nabla \Psi(x)|^2 |\eta_\epsilon(x)|^2 \mathrm{d} x + \frac{1}{\epsilon^2} \int_{V_\epsilon \cap I_N} |\Psi(x)|^2 \mathrm{d} x,
    \end{align*}
    where $V_\epsilon \coloneqq \{x: d_V(x) \leq 2\epsilon\}$. The first term converges to zero again by dominated convergence. For the second term, we can assume, via a partition of the unity argument, that $V_\epsilon \cap I_N$ is the hypergraph of a Lipschitz function, i.e., up to a rotation and translation, there exists $\kappa:\R^{N-1}\rightarrow \R$ Lipschitz such that $\partial I_N \cap V_\epsilon = \{(x',\kappa(x')): x' \in W_\epsilon\}$ for some open bounded set $W_\epsilon \subset \R^{N-1}$ and $V_\epsilon \cap I_N = \{(x',x_N) : x' \in W_\epsilon, \quad  \kappa(x') < x_N < y(x',\epsilon)\}$, where $y(x',\epsilon)-\kappa(x')\lesssim \epsilon$. Since $V_\epsilon \cap \partial I_N \subset U$ for $\epsilon>0$ small enough, we have $\Psi(x',\kappa(x')) = 0$, and therefore
    \begin{align*}
        \frac{1}{\epsilon^2} \int_{V_\epsilon} |\Psi(x)|^2 \mathrm{d} x = \frac{1}{\epsilon^2} \int_{W_\epsilon} \mathrm{d} x' \int_{\kappa(x')}^{y(x',\epsilon)} |\Psi(x',x_N)|^2 \mathrm{d} x_N = \frac{1}{\epsilon^2} \int_{W_\epsilon} \int_{\kappa(x')}^{y(x',\epsilon)} \left|\int_{\kappa(x')}^{x_N} \partial_{x_N} \Psi(x',t)\mathrm{d} t\right|^2 \mathrm{d} x.
    \end{align*}
    If we now apply Cauchy-Schwarz, use the change of variables $z = (x_N-\kappa(x'))/\epsilon$, and use the estimate $y(x',\epsilon)-\kappa(x') \leq C\epsilon$, we find
    \begin{align*}
        \frac{1}{\epsilon^2} \int_{W_\epsilon} \int_{\kappa(x')}^{y(x',\epsilon)} \left|\int_{\kappa(x')}^{x_N} \partial_{x_N} \Psi(x',t) \mathrm{d} t\right|^2 \mathrm{d} x&\leq \int_{W_\epsilon} \int_{\kappa(x')}^{y(x',\epsilon)} \left(\int_{\kappa(x')}^{x_N} |\partial_{x_N} \Psi(x',t)|^2 \mathrm{d} t \right) \frac{x_N-\kappa(x')}{\epsilon^2} \mathrm{d} x_N \mathrm{d} x' \\
        &\lesssim \int_{W_\epsilon} \int_0^{C} \int_{\kappa(x')}^{\kappa(x')+\epsilon z} |\partial_{x_N} \Psi(x',t)|^2 \mathrm{d} t \mathrm{d} z \mathrm{d} x'\\
        &\lesssim \int_{V_\epsilon\cap I_N} |\partial_{x_N} \Psi(x',t)|^2 \mathrm{d} x'\mathrm{d} t.
    \end{align*}
    Hence, by dominated convergence, the limit $\epsilon \downarrow 0^+$ goes to zero, which completes the proof.
\end{proof}

We now show that the Neumann trace is meaningful along any open subset of the boundary where $\Psi$ vanishes. The key to prove this is the locality of the form $\mathfrak{a}_{v,w}$.
\begin{lemma}[Neumann trace formula] \label{lem:Ntrace formula} Let $U\subset \partial I_N$ be a relatively open subset of $I_N$ and suppose that $U$ is compactly contained in the interior of a face of $\partial I_N$. Let $\Psi$ be an eigenfunction of $H_N^\alpha(v,w;\Gamma)$ such that $(\gamma\Psi) \rvert_U = 0$. Then for any $f \in C^\infty_c(U)$, we have
\begin{align}
(\partial_\nu \Psi)(f) = -\lim_{\epsilon \downarrow 0^+} \frac{\inner{\gamma_\epsilon \Psi,f}_{L^2(U)}}{\epsilon} , \label{eq:trace formula}
\end{align}
where the above limit exists and $\gamma_\epsilon \Psi$ denotes the Dirichlet trace of $\Psi$ along the hyperplane parallel and at distance $\epsilon$ of the face containing $U$.
\end{lemma}

\begin{proof} Without loss of generality, let us assume that $U\subset \subset E_0 \coloneqq \{ x =(x_1,x')\in \partial I_N : x_1 = 0, \quad x' \in I_{N-1}\}$. Then, we let $\beta \in C^\infty_c(\R;\R)$ be a standard mollifier, i.e., $0\leq \beta \leq 1$, $\mathrm{supp}(\beta) \subset [-1,1]$ and $\beta = 1$ on a neighborhood of $0$, and define the function $F$ as
    \begin{align*}
        F(x_1,x') \coloneqq \beta(x_1) f(x')\quad \mbox{for $(x_1,x') \in I_N$.}
    \end{align*}
    Clearly, $F$ is an $\mH^1$-extension of $f$, i.e., $F\in \mH^1(I_N)$ and $\gamma F = f$.
    
    Next, let $V\subset \partial I_N$ be such that $\mathrm{supp}(f) \subset \subset V \subset \subset U$ and define $\eta_\epsilon$ as in Lemma~\ref{lem:H10 approx}. Then $F_\epsilon \coloneqq (1-\eta_\epsilon) F \in \mH^1_0(I_N)$ and therefore
    \begin{align}
        0 = \mathfrak{a}_{v,w}(\Psi, P_N F_\epsilon)- \lambda \inner{\Psi, P_N F_\epsilon}  =&\int_{I_N} \left(1-\eta_\epsilon\right) \overline{\nabla \Psi} \scpr \nabla F  \mathrm{d} x +  v(\rho_{\Psi, P_N F_\epsilon}) + w(\rho^{(2)}_{\Psi, P_N F_\epsilon}) - \lambda\inner{\Psi,F_\epsilon} \label{eq:middleeq} \\
        &- \int_{I_N} F(x) \overline{\nabla \Psi(x)} \scpr \nabla \eta_\epsilon(x) \mathrm{d} x, \label{eq:lasteq}
    \end{align}
    where the first equality follows because $\Psi$ is an eigenfunction of $H_N^\alpha(v,w;\Gamma)$ with eigenvalue $\lambda$.
    
    The key observation now is the following: since
    \begin{align*}
        \overline{\Psi(x)} (P_NF_\epsilon)(x) = \frac{1}{N!} \sum_{\sigma \in \mathcal{P}_N} \mathrm{sgn}(\sigma) \overline{\Psi(x) \left(1-\eta_\epsilon(\sigma x)\right)} F(\sigma x),
    \end{align*}
    and $\gamma\Psi$ vanishes on $\sigma^{-1}(U)$ for any $\sigma \in \mathcal{P}_N$ (by antisymmetry), from Lemma~\ref{lem:H10 approx} we have
    \begin{align*}
        \lim_{\epsilon \downarrow 0^+} \norm{\Psi(1-\sigma_\#\eta_\epsilon) - \Psi}_{\mH^1} = 0 \quad \mbox{for any $\sigma \in \mathcal{P}_N$.}
    \end{align*}
    Consequently, by Lemma~\ref{lem:regularity reduced densities}, 
    \begin{align*}
        \lim_{\epsilon \downarrow 0^+} \rho_{\Psi,P_NF_\epsilon} = \rho_{\Psi, P_N F} \quad \mbox{in $\mH^1(I)$} \quad \mbox{and}\quad  \lim_{\epsilon \downarrow 0^+} \rho^{(2)}_{\Psi, P_N F_\epsilon} = \rho^{(2)}_{\Psi, P_N F} \quad \mbox{in $\mW^{-1,p}(I_2)$.}
\end{align*}
Using these convergence results and dominated convergence for the kinetic term, we see that equation~\eqref{eq:middleeq} converges to $\partial_\nu \Psi(f)$ as $\epsilon\downarrow 0^+$. In particular, 
    \begin{align*}
        \partial_\nu \Psi(f) = \lim_{\epsilon \downarrow 0^+}  \int_{I_N} F(x) \overline{\nabla \Psi(x)} \scpr \nabla \eta_\epsilon(x) \mathrm{d} x,
    \end{align*}
where the limit exists. To compute this limit, let us take
\begin{align*}
    \eta(x) = \begin{dcases} 1, \quad \mbox{for $x\leq 1/2$,}\\
    2-2x, \quad \mbox{for $\frac12 \leq x \leq 1$,}\\
    0, \quad \mbox{otherwise.} 
    \end{dcases}
\end{align*}
Since $\mathrm{supp}(F) \subset [0,1] \times V$, we have $d_V(x)= x_1$ for any $x =(x_1,x') \in \mathrm{supp}(F)$. In particular, $\eta_\epsilon(x) = \eta(x_1/\epsilon)$ for $x\in \mathrm{supp}(F)$. Moreover, for $\epsilon>0$ small enough, we have $\beta(x_1) = 1$ on $\mathrm{supp}(\eta_\epsilon(x))$. Therefore
\begin{align}
    \partial_\nu \Psi(f) &= \lim_{\epsilon\downarrow 0^+} \int_{I_N} F(x) \overline{\nabla \Psi(x)} \scpr \nabla \eta_\epsilon(x) \mathrm{d} x \nonumber \\
    &= \lim_{\epsilon \downarrow 0^+} \int_{ [\epsilon/2,\epsilon] \times V} f(x')  \overline{\partial_{x_1} \Psi(x_1,x')}  \left(\frac{-2}{\epsilon}\right)\mathrm{d} x_1 \mathrm{d} x'\nonumber \\
    &= \lim_{\epsilon \downarrow 0^+} \frac{\inner{\gamma_{\epsilon/2}\Psi,f} - \inner{\gamma_{\epsilon} \Psi,f}}{\epsilon/2} \label{eq:limit}
\end{align}
To complete the proof, we now write
\begin{align*}
    \inner{\gamma_{2^{-n}\epsilon}\Psi,f} - \inner{\gamma_\epsilon \Psi, f}  - (\epsilon-2^{-n} \epsilon) \partial_\nu \Psi(f) = \sum_{k=1}^{n} \left(\inner{\gamma_{2^{-k} \epsilon} \Psi,f} - \inner{\gamma_{2^{-k+1}\epsilon} f} - 2^{-k}\epsilon \partial_\nu \Psi\right)
\end{align*}
Thus, by~\eqref{eq:limit}, for any $\delta>0$ we can take $\epsilon>0$ so small that the right hand side is controlled by $\delta \epsilon$. Taking the limit $n \ra \infty$ and using that $\gamma \Psi = 0$ on $U$, we conclude that
\begin{align*}
    |\epsilon \partial_\nu \Psi(f) + \inner{\gamma_\epsilon \Psi,f}| \lesssim \delta \epsilon. 
\end{align*}
Diving this expression by $\epsilon$ and taking the limits $\epsilon \downarrow 0^+$ and $\delta \downarrow 0^+$ completes the proof.
\end{proof}

We now show that the ground-state $\Psi$ of $H_N^\alpha(v,w;\Gamma)$ satisfies $\partial_\nu \Psi = 0$ along any open set $U\subset \partial S_N$ such that $\gamma \Psi \rvert_U = 0$. To prove this, the key observation is that $\Psi$ is non-negative in $S_N$ by the Perron-Frobenius theorem. 
\begin{lemma}[Vanishing of normal derivative] \label{lem:trace zero} Let $\alpha \in \R\setminus \{0\}$ satisfy~\eqref{eq:parity condition} and let $\Psi$ be the ground-state of $H_N^{\alpha}(v,w;\Gamma)$. Suppose that $(\gamma \Psi)\rvert_U = 0$ for some (relatively) open set $U\subset \subset \mathrm{int}\, F_0 \setminus \Gamma_N$, where $F_0 =\{ (x_1,x_2,...,x_N) \in \partial S_N : x_1 = 0\}$. Then 
\begin{align}
    \partial_\nu \Psi(f)= 0\quad \mbox{for any $f\in \mH^{1/2}(\partial \Omega)$ with $\mathrm{supp}(f) \subset U$.} \label{eq:zero trace}
\end{align}
\end{lemma}

\begin{proof} Since any function $f\in \mH^{1/2}(\partial I_N)$ with $\mathrm{supp}(f) \subset U$ can be approximated by functions in $C_c^\infty(U)$ it suffices to prove~\eqref{eq:zero trace} for $f\in C^\infty_c(U)$. To this end, we first claim that
\begin{align}
    \partial_\nu \Psi(f) \leq 0.  \quad \mbox{for any $f\in C^\infty_c(U)$ such that $f \geq 0$.} \label{eq:negativity}
\end{align}
Indeed, since $\Psi \geq 0$ almost everywhere in $S_N$ (by Theorem~\ref{thm:PFFermi}), this claim follows from the formula in Lemma~\ref{lem:Ntrace formula}.

Next, we recall that, by the definition of $\mH^1_\alpha(I_N)$ and antisymmetry, any $\Psi \in \mathcal{Q}_N^\alpha$ satisfies 
\begin{align*}
    (\gamma \Psi)(0,x') - (-1)^{N-1} \alpha (\gamma \Psi)(x',1) = 0 \quad\mbox{for a.e. $x' \in I_{N-1}$.}
\end{align*} 
Therefore $\gamma \Psi$ also vanishes on the set $U' \coloneqq \{(x',1) : (0,x') \in U\} \subset F_N$. Consequently,~\eqref{eq:negativity} also holds for $f\in \mH^{1/2}(\partial I_N)$ with $\mathrm{supp}(f) \subset U'$. Now let $f\in C^\infty_c(U)$ and define
\begin{align*}
    f'(x) \coloneqq \begin{cases} (-1)^{N-1}\alpha f(0,x') \quad &\mbox{for $x=(x',1) \in U'$,}\\
    0, \quad &\mbox{for $x\in \partial I_N \setminus U'$.} \end{cases}
\end{align*}
Let $Jf$ and $Jf'$ be extensions of $f$ and $f'$ in $\mH^1(I_N)$. Then, recalling the proof of Lemma~\ref{lem:positivespaces}, we see that $P_N(Jf+Jf') \in \mathcal{Q}_N^\alpha(\Gamma)$. Therefore, by the definition of the Neumann trace, we have
\begin{align*}
    \partial_\nu \Psi(f) + \partial_\nu \Psi(f') &= \partial_\nu \Psi(f+f') = \mathfrak{a}_{v,w}\left(\Psi,P_N(Jf+Jf')\right) - \lambda \inner{\Psi, P_N(Jf+Jf')} = 0 .
\end{align*}
Thus if $f\geq 0$, then $f'\geq 0$ by~\eqref{eq:parity condition}, and therefore $\partial_\nu \Psi(f) = \partial_\nu \Psi(f') = 0$ by~\eqref{eq:negativity}. This concludes the proof.
\end{proof}

We can now complete the proof of Theorem~\ref{thm:UCP boundary}.
\begin{proof}[Proof of Theorem~\ref{thm:UCP boundary}: the case of non-local BCs] Let $U\subset \partial I_N \setminus \Gamma_N$ be relatively open. We shall assume that $\Psi \rvert_U = 0$ and obtain a contradiction. To this end, first note that, by the antisymmetry of $\Psi$, the Dirichlet trace $\gamma \Psi$ vanishes on the symmetrized set $U_N = \cup_{\sigma} \sigma(U)$. Moreover, after possibly shrinking $U_N$, we can assume that $U_N \cap \partial S_N$ is compactly contained on the interior of an exterior face of $\partial S_N$. So without loss of generality, we assume that $U \subset  F_0 =\{(x_1,...,x_N)\in \partial S_N : x_1 =0\}$. Hence, if we denote by $\widetilde{\Psi} = T^{-1} \Psi = \sqrt{N!} \Psi\rvert_{S_N}$ the ground-state of the reduced operator $\widetilde{H}_N^\alpha(v,w;\Gamma)$, Lemma~\ref{lem:trace zero} implies that
\begin{align*}
    \partial_\nu \widetilde{\Psi}(f) &\coloneqq \widetilde{\mathfrak{a}}_{v,w}(\widetilde{\Psi},Jf) - \lambda \inner{\widetilde{\Psi},Jf} = \mathfrak{a}_{v,w}(T\widetilde{\Psi}, T Jf) - \lambda \inner{T \widetilde{\Psi}, TJf} = \sqrt{N!} \partial_\nu \Psi(f) = 0,
\end{align*}
for any $f \in \mH^{1/2}(\partial S_N)$ such that $\mathrm{supp}(f) \subset U$, where $J: \mH^{1/2}(\partial S_N) \rightarrow \mH^1(S_N)$ is any right inverse of the Dirichlet boundary trace on $S_N$.

We now claim that this is equivalent to
\begin{align}
    \widetilde{\mathfrak{a}}_{v,w}(\widetilde{\Psi}, \Phi) - \lambda\inner{\widetilde{\Psi}, \Phi} = 0, \quad \mbox{for any $\Phi \in \mathcal{Q}_W$,} \label{eq:claim2}
\end{align}
where
\begin{align*}
    \mathcal{Q}_W \coloneqq \{ \Phi \in \mathcal{Q}_N(\Gamma;S_N):  (\gamma \Phi)(0,x') - (-1)^{N-1}\alpha (\gamma \Phi)(x',1)\quad \mbox{for $(0,x') \in F_0 \setminus W$} \},
\end{align*}
for any neighborhood $W\subset U \subset F_0$. Indeed, if $\Phi \in \mathcal{Q}_W$, then $\Phi$ can be written as the sum of a function in the form domain of the reduced operator $\widetilde{H}^\alpha_N(v,w;\Gamma)$
and a function with trace supported in $\overline{W} \subset U$. To be precise, for any $\Phi \in \mathcal{Q}_W$, we can define
\begin{align*}
    f_{\Phi} (x) \coloneqq \begin{dcases} (-1)^{N-1}\frac{1}{\alpha} (\gamma \Phi)(x',1), \quad &\mbox{if $x=(0,x')$ for $x' \in S_{N-1}$},\\
    (\gamma \Phi)(x',1),\quad &\mbox{if $x=(x',1)$ with $x'\in S_{N-1}$,}\\
    0,\quad &\mbox{otherwise,}\end{dcases}
\end{align*}
and note that $f_{\Phi}\in \mH^{1/2}(\partial S_N)$ because $f_{\Phi}$ is locally equal to $\gamma \Phi$ or a rigid motion of $\gamma \Phi$. Moreover, by the definition of $\mathcal{Q}_W$, we have $\mathrm{supp}(\gamma \Phi - f_\Phi) \subset \overline{W}$. Thus, since $J f_\Phi$ belongs to the form domain of $\widetilde{H}_N^\alpha(v,w;\Gamma)$ (namely~\eqref{eq:claim} intersected with $\mH^1(S_N;\Gamma_N\cap \partial S_N)$), we have the desired decomposition $\Phi = (\Phi - Jf_\Phi) + Jf_\Phi$, which proves the claim in~\eqref{eq:claim2}. 

Next, notice that~\eqref{eq:claim2} is equivalent to saying that $\widetilde{\Psi}$ is an eigenfunction of the operator associated to the form $\widetilde{\mathfrak{a}}_{v,w}$ but with the larger form domain $\mathcal{Q}_W$. The important observation now is that the trace of functions in $\mathcal{Q}_W$ have no restriction inside $W$. Therefore, we can repeat the extension argument in the proof of the case of local BCs. More precisely, we extend $S_N$ to $\Omega$ by adding an open ball whose intersection with $\partial S_N$ is compactly contained in $W$. Then $\widetilde{\Psi}$ extended by zero is also an eigenfunction of the extended operator $\widetilde{H}_\Omega$ associated to~\eqref{eq:form Omega} but this time with form domain 
\begin{align*}
    \mathcal{Q}_\Omega \coloneqq \{ \Psi \in \mH^1(\Omega;\Gamma) : \Psi\rvert_{U\cap \partial \Omega} = 0 \quad \mbox{and}\quad (\gamma \Phi)(0,x') - (-1)^{N-1}\alpha(\gamma \Phi)(x',1)\quad \mbox{for $(0,x') \in F_0 \setminus W$} \}.
\end{align*}
Thus, recalling assumption~\eqref{eq:parity condition}, one can repeat the arguments in the proof of Lemma~\ref{lem:positivespaces} to show that $\mathcal{Q}_\Omega$ is positivity and reality preserving. In particular, Theorem~\ref{thm:PF} applies to $\widetilde{H}_\Omega$ and therefore its ground-state $\Psi_\Omega$ is non-degenerate and almost everywhere strictly positive in $\Omega$. This in turn implies that $\Psi_\Omega$ is either orthogonal or parallel to $\widetilde{\Psi}$, which gives us a contradiction as these functions can not be orthogonal since they are strictly positive in $S_N$ but can not be parallel because $\widetilde{\Psi}$ vanishes identically in $\Omega \setminus S_N$ while $\Psi_\Omega$ does not.
\end{proof}

\subsection{Proof of Theorems~\ref{thm:eigenvalues}}

We can now complete the proof of Theorem~\ref{thm:eigenvalues}.

\begin{proof}[Proof of Theorem~\ref{thm:eigenvalues}]
As mentioned in Remark~\ref{rem:strict}, the inequality $\lambda_1(\Gamma)\leq \lambda_1(\Gamma')$ is immediate from the variational principle. Hence, it suffices to show that equality cannot hold. For this, note that, if $\lambda_1(\Gamma) = \lambda_1(\Gamma')$, then by the variational principle and the non-degeneracy theorem, the ground-state of $H^{\alpha}_N(v,w;\Gamma)$ and $H_N^\alpha(v,w;\Gamma')$ are the same. Hence the ground-state $\Psi$ of $H_N^\alpha(v,w;\Gamma)$ vanishes on $\Gamma'$. As $\Psi$ is antisymmetric, it must also vanish in the symmetrized set $\Gamma_N'$. Hence, $\Psi$ vanishes on $\Gamma_N'\setminus \Gamma_N$. As this set contains a relatively open subset of the boundary by assumption, this is not possible by Theorem~\ref{thm:UCP boundary}. We thus have a contradiction, which concludes the proof.
\end{proof}

\addtocontents{toc}{\protect\setcounter{tocdepth}{1}}
\addtocontents{toc}{\protect\setcounter{tocdepth}{-1}}
\section*{Acknowledgements}
The author wishes to thank Chokri Manai for several inspiring discussions, which led to the ideas in the proof of Theorem~\ref{thm:UCP boundary}. The author also thanks Timo Weidl, Robert van Leeuwen, and Markus Penz for fruitful discussions and helpful comments on this article. The author is also grateful to Rafael Antonio Lainez Reyes for pointing out the references on Sobolev spaces. In addition, the author is grateful to the anonymous referees for several suggestions to improve the presentation of the paper.

T.C.~Corso acknowledges funding by the \emph{Deutsche Forschungsgemeinschaft} (DFG, German Research Foundation) - Project number 442047500 through the Collaborative Research Center "Sparsity and Singular Structures" (SFB 1481). 

\addtocontents{toc}{\protect\setcounter{tocdepth}{1}}


\section*{Data availability}
No datasets were generated or analysed during the current study.

\section*{Competing interests}

The authors has no competing interests to declare that are relevant to the content of this article.


\appendix

\section{Schr\"odinger operators with distributional potentials}

In this section, we present the construction of self-adjoint realizations of $H_N(v,w)$ via their associated quadratic form. We also list some important examples of distributional potentials for which this construction holds. A detailed exposition on the theory quadratic forms and self-adjoint extensions can be found in several classical references \cite{RS75,RS78,AGHH88,Tes14}.

\subsection{Self-adjoint realization} \label{sec:self-adjoint} The estimate in Lemma~\ref{lem:regularity reduced densities} allows us to define self-adjoint Schr\"odinger operators with external and interaction potentials in a large class of distributions. Indeed, from estimate~\eqref{eq:reduced Sobolev} and Young's inequality, we find that any $v\in \mathcal{V}$ and $w\in \mathcal{W}$ satisfy
\begin{align}
    |v(\rho_{\Psi,\Psi}) + w(\rho^{(2)}_{\Psi,\Psi})| \leq \epsilon  \norm{\Psi}_{\mH^1}^2 + C_\epsilon \norm{\Psi}_{\mL^2}^2, \quad \mbox{for any $\Psi \in \mH^1(I_N)$,}  \label{eq:KLMNest}
\end{align}
and any $\epsilon>0$ for some constant $C_\epsilon = C(v,w,\epsilon)>0$. Hence, from the celebrated KLMN theorem (see \cite[Theorem X.17]{RS75}), for any closed subspace $B \subset \mH^{1/2}(\partial \Omega)$, we can construct a unique self-adjoint operator associated to the sesquilinear form
\begin{align}
    \mathfrak{a}_{v,w}(\Psi,\Phi) = \int_{I_N} \overline{\nabla \Psi(x)} \scpr \nabla \Phi(x) \mathrm{d} x + v(\rho_{\Psi,\Phi}) + w(\rho^{(2)}_{\Psi,\Phi}), \quad \mbox{for $\Psi,\Phi \in \mH^1_B(I_N)\cap \mathcal{H}_N$.} \label{eq:quadratic form}
\end{align}
Precisely, the following holds.
\begin{lemma}[Schr\"odinger operators with distributional potentials] \label{lem:operator def} Let $N\in \N$, $v\in \mathcal{V}$, and $w\in \mathcal{W}$. Moreover, let $B\subset \mH^{1/2}(\partial I_N)$ be a closed subspace and $\mH^1_B(I_N)$ the trace-restricted Sobolev space introduced in Definition~\ref{def:traceSobolev}. Then the quadratic form in~\eqref{eq:quadratic form} is closed, symmetric and semibounded. In particular, there exists a unique self-adjoint operator $H_N^B(v,w)$ associated to this form. Moreover, $H_N^B(v,w)$ is semibounded and has purely discrete spectrum.
\end{lemma}

\begin{remark}[Discrete spectrum via compact form domain]\label{rem:discrete spectrum} Recall that $H_N^B(v,w)$ has discrete specturm if and only if the form domain is compactly embedded in $\mathcal{H}_N$. Hence, the discrete spectrum statement in Lemma~\ref{lem:operator def} follows from the compact embedding $\mH^1_B(I_N) \subset \mH^1(I_N) \hookrightarrow \mL^2(I_N)$.
\end{remark}

\begin{remark}[Notation]
Throughout the paper, we denote by $H_N(v,w;\Gamma)$ the self-adjoint realization of $\mathfrak{a}_{v,w}$ with form domain $\mH^1(I_N;\Gamma) \cap \mathcal{H}_N$ and by $H^\alpha_N(v;w;\Gamma)$ the realization with form domain $\mH^1_\alpha(I_N;\Gamma)\cap \mathcal{H}_N$. 
\end{remark}

\subsection{From single-particle to many-particle boundary conditions} Formally, the non-interacting operator $H_N(v,0)$ can be written as the sum 
\begin{align*}
    H_N(v,0) = \sum_{j=1}^N 1\otimes... \otimes \overbrace{h(v)}^{\mathclap{j^{th} position}} \otimes ... \otimes 1, 
\end{align*}
where $1$ is the identity operator and $h(v)= -\Delta +v$ is the single-particle operator acting on $\mL^2(I)$. This allow us to reduce the study of non-interacting (electronic) operators to the case of single-particle ones. However, to make this observation precise, we need to impose appropriate boundary conditions on the form domains of $H_N(v,0)$ and $h(v)$. In this section, our goal is to clarify this point. 

To this end, we start with the following lemma.
\begin{lemma}[Many-particle form domain]\label{lem:density single-to-many} Let $\alpha \in \R$ and $\mH^1_\alpha(I) = \{ f\in \mH^1(I): f(0) =\alpha f(1) \}$. Then for any $N\in \N$, the span of Slater determinants
\begin{align*}
    \Phi = \phi_1 \wedge... \wedge \phi_N, \quad \phi_j \in \mH^1_\alpha(I)
\end{align*}
is dense in $\mathcal{Q}_N^\alpha = \mH^1_\alpha(I_N)\cap \mathcal{H}_N$.
\end{lemma}

\begin{proof} We first observe that the projection on the antisymmetric space, explicitly given in~\eqref{eq:antisymmetric projection}, is continuous with respect to the $\mH^1$ topology and maps $N$-fold tensor products to $N$-fold Slater determinants. Therefore, it suffices to show that the space
\begin{align*}
    Y \coloneqq \mathrm{span}\{ \phi_1 \otimes... \otimes \phi_N : \phi_j \in \mH^1_\alpha(I) \}
\end{align*}
is dense in $\mH^1_\alpha(I_N)$. For this, let $\{\phi_j\}_{j \in \N}$ be an orthonormal basis of eigenfunctions of the one-dimensional Laplacian in $\mH^1_\alpha(I)$, i.e., $\phi_j$ satisfies
\begin{align*}
    \int_I \overline{\partial_x \phi_j(x)} \partial_x \psi(x) \mathrm{d} x  = \lambda_j \int_I \overline{\phi_j(x)} \psi(x) \mathrm{d} x , \quad \mbox{for any $\psi \in \mH^1_\alpha(I)$.}
\end{align*}
Then, it is not hard to see that $\phi_j \in C^\infty(\overline{I})$ and 
\begin{align*}
    \alpha \partial_x \phi_j(0) -\partial_x \phi_j(1) = 0. 
\end{align*}
Consequently, by integration by parts, for any $\Psi \in \mH^1_\alpha(I_N)$ and $(j_1,...,j_N) \in \N^N$ we have
\begin{align*}
    \inner{\Psi, \phi_{j_1}\otimes... \otimes \phi_{j_N}}_{\mH^1} &= \sum_{k=1}^N \inner{\partial_{x_k} \Psi, \partial_{x_k} \phi_{j_1} \otimes... \phi_{j_N}}_{\mL^2} + \inner{\Psi, \phi_{j_1}\otimes ... \phi_{j_N}}_{\mL^2}\\
    &=  \left(1+\sum_{k=1}^N \lambda_{j_k}\right) \inner{\Psi, \phi_{j_1} \otimes... \otimes \phi_{j_N}}.
\end{align*}
Thus, if $\Psi$ is $\mH^1$-orthogonal to $Y$, then $\Psi$ is $\mL^2$-orthogonal to all $\{\phi_{j_1} \otimes... \otimes \phi_{j_N}\}_{(j_1,...,j_N) \in \N^N}$. As this set is an orthogonal basis of $\mL^2(I_N)$, we must have $\Psi = 0$. Hence, the intersection of $\mH^1_\alpha(I_N)$ with the $\mH^1$-orthogonal complement of $Y$ is trivial, which implies that $Y$ is dense in $\mH^1_\alpha(I_N)$. 
\end{proof}

As a consequence of the previous lemma, we have the following relation between the single-particle operator $h_\alpha(v)\coloneqq H_1^\alpha(v,0)$ and the $N$-particle operator $H_N^\alpha(v,0)$ associated to the form
\begin{align*}
    \mathfrak{a}_{v}(\Psi,\Psi) = \int_{I_N} |\nabla \Psi(x)|^2 \mathrm{d} x + v(\rho_{\Psi}), \quad \Psi \in \mH^1_\alpha(I_N) \cap \mathcal{H}_N.
\end{align*}
\begin{lemma}[From single-particle to non-interacting many-particle systems] \label{lem:single-to-many} Let $v\in \mathcal{V}$, $\alpha \in \C$, and $N \in \N$. Let $\{\psi_j\}_{j \in \N}$ be an orthogonal basis of eigenfunctions\footnote{Note that such a basis exists because the spectrum of $h_\alpha(v)$ is purely discrete, see Remark~\ref{rem:discrete spectrum}.} of $h_\alpha(v) = H_1^\alpha(v,0)$ with corresponding eigenvalues $\{\lambda_j\}_{j\in \N}$ ordered in non-decreasing order. Then the Slater determinants
\begin{align}
    \{ \Psi_J = \psi_{j_1} \wedge... \wedge \psi_{j_N}  : J = \{j_1,...,j_N\} \in \N^N \quad \mbox{such that} \quad 1\leq j_1 < j_2 ... <j_N\} \label{eq:Slater}
\end{align}
forms an $\mathcal{H}_N$-orthogonal basis of eigenfunctions of $H^\alpha_N(v,0)$ with eigenvalues $\lambda_J = \sum_{k=1}^N \lambda_{j_k}$. In particular, the ground-state of $H_N^\alpha(v,0)$ is non-degenerate if and only if $\lambda_N$ satisfies the closed shell condition $\lambda_N < \lambda_{N+1}$.
\end{lemma}

\begin{proof} From a long but straightforward calculation we find that
\begin{align*}
    \mathfrak{a}_v(\Psi_J, \phi_1 \wedge... \wedge \phi_N)= \lambda_J \inner{\Psi_J, \phi_1 \wedge... \wedge \phi_N}_{\mL^2}, \quad \mbox{for any $\{\phi_j\}_{j=1}^N \subset \mH^1_\alpha(I)$,}
\end{align*}
where  $\Psi_J = \psi_{j_1} \wedge... \wedge \psi_{j_N}$ and $\lambda_J = \sum_{k=1}^N \lambda_{j_k}$. Thus by Lemma~\ref{lem:density single-to-many} and approximation we have
\begin{align*}
    \mathfrak{a}_v(\Psi_J, \Phi) = \lambda_J \inner{\Psi_J,\Phi} \quad \mbox{for any $\Phi \in \mH^1_\alpha(I_N)\cap \mathcal{H}_N$.}
\end{align*}
Therefore $\Psi_J$ is an eigenfunction of $H_N^\alpha(v,0)$ with eigenvalue $\lambda_J$. The fact that the $\Psi_J$'s form an $\mathcal{H}_N$-orthogonal basis follows from the fact that $\{\psi_j\}_{j \in \N}$ is an $\mL^2$-orthogonal basis of $\mL^2(I)$. 
\end{proof}

\subsection{Examples of distributional potentials}\label{sec:examples} To illustrate how large the class of admissible external and interaction potentials in Lemma~\ref{lem:operator def} is, let us present a few important examples.

\begin{enumerate}[label=(\arabic*)]
    \item The Dirac's delta potential $\delta_{x_0}$ with $x_0 \in I$ belongs to $\mathcal{V}$ by the continuous embedding $\mH^1(I) \subset C(\overline{I})$. Note that $x_0 \in \partial I = \{0,1\}$ is also allowed. 
    \item The $\delta$-interaction potential $w = \delta(x-y)$, defined via
    \begin{align*}
        w_\delta(\rho^{(2)}_{\Psi,\Phi}) = \int_I \rho^{(2)}_{\Psi,\Phi}(x,x) \mathrm{d} x.
    \end{align*}
    This follows from the fact that the trace operator $\gamma: \mW^{1,p}(I_2) \rightarrow \mL^p(D)$ sending $\rho_2$ to its restriction along the diagonal set $D =\{(x,x) \in I_2:x\in I\}$ is a continuous operator. 
    \item Multiplicative potentials $v\in \mL^1(I)$ and $w\in \mL^1(2I)$, whose actions are defined via
    \begin{align*}
        v(\rho) = \int_I v(x) \rho(x) \mathrm{d} x \quad \mbox{and}\quad w(\rho_2) = \int_{I_2} w(x-y) \rho_2(x,y) \mathrm{d} x \mathrm{d}y,
    \end{align*}
    are also allowed. The first one is clear from H\"older's inequality. For the second one, note that
    \begin{align*}
        f(u) \coloneqq \int_{|u|}^{1-|u|}\rho_2\left(\frac{u+v}{2},\frac{u-v}{2}\right)\mathrm{d} v 
    \end{align*}
    belongs to $\mL^\infty(I)$ for any $\rho_2 \in \mW^{1,p}(I_2)$ with $1\leq p \leq \infty$. Hence, the map
    \begin{align*}
        \rho_2 \mapsto \int_{I_2} w(x-y) \rho_2(x,y) \mathrm{d}x \mathrm{d} y = \frac{1}{2} \int_{-1}^1 w(u) f(u) \mathrm{d} u
    \end{align*}
    belongs to $\mW^{-1,q}(I_2)$. 
    \item Lemma~\ref{lem:operator def} can also be extended to the case of $k^{th}$-body distributional potentials for $k \geq 3$. To be precise, for any $w_k\in \mW^{-1,q}(I_k)$ with $q>k$, one can use Lemma~\ref{lem:regularity reduced densities} and Young's inequality to show that the form $(\Psi,\Phi) \mapsto w_k(\rho^{(k)}_{\Psi,\Phi})$, where
    \begin{align*}
        \rho^{(k)}_{\Psi,\Phi}(x_1,..,x_k) = \frac{N!}{(N-k)! k!} \int_{I_{N-k}} \overline{\Psi(x_1,...,x_N)} \Phi(x_1,...,x_N)\mathrm{d}x_{k+1} ... \mathrm{d} x_N,
    \end{align*}
    is $\Delta$-bounded with relative bound $0$. Therefore, any Schr\"odinger operator of the form
    \begin{align*}
       H_N(w_1,w_2,...,w_N) = -\Delta + \sum_{k=1}^N \sum_{j_1\neq ... \neq j_k}^N w_k(x_{j_1},x_{j_2},...,x_{j_k}), 
    \end{align*}
    with real-valued $w_k \in \mW^{-1,q_k}(I_k)$, where $q_k >k$ for $k\geq 2$ and $w_1 \in \mathcal{V}$, defines a self-adjoint operator with discrete spectrum. In particular, the case of 3D Coulomb interactions
    \begin{align*}
        w(x_1,...,x_6) = \left(\sum_{j=1}^3 (x_j-x_{j+3})^2\right)^{-\frac12}
    \end{align*}
    is also included for $N \geq 6$.
    \end{enumerate}

\section{A prototypical example: the one-dimensional Laplacian}
\label{app:Laplace}
In this section we gather some well-known formulas for the eigenfunctions of the one-dimensional Laplacian under different boundary conditions. These examples serve to illustrate the results in Theorems~\ref{thm:main},~\ref{thm:non-local}, and~\ref{thm:non-interacting} and to justify Remark~\ref{rem:parity condition}.

Throughout this section we denote respectively by $-\Delta$, $-\Delta_0$, and $-\Delta_\pm$ the self-adjoint realizations of the one-dimensional Laplacian on $I$ with Neumann, Dirichlet, and periodic/anti-periodic boundary conditions. Precisely, these are the unique self-adjoint operators associated to the form
\begin{align*}
    \mathfrak{a}_0(\psi,\psi) = \norm{\nabla \psi}_{\mL^2}^2,
\end{align*}
with form domains $\mathcal{Q} = \mH^1(I)$, $\mathcal{Q}_0 = \mH^1_0(I)$, and $\mathcal{Q}_\pm = \mH^1_{\pm 1}(I)$. The eigenvalues and eigenfunctions of these operators can be explicitly computed.
\begin{proposition}[Spectrum of one-dimensional Laplacian] \label{prop:1D Laplacian} The following holds: \begin{enumerate}[label=(\roman*)]
\item(Neumann) We have
\begin{align*}
    \sigma(-\Delta) = \{  \pi^2 (k-1)^2  : k\in \N\}.
\end{align*}
Moreover, each eigenvalue is simple and the associated normalized eigenfunction are
\begin{align*}
    \phi_k(x) = \sqrt{2} \cos(\pi (k-1) x), \quad k \in \N \quad \mbox{and}\quad \phi_1(x) = 1.
\end{align*}
\item(Dirichlet) We have
\begin{align*}
    \sigma(-\Delta_0) = \{ \pi^2 k^2: k \in \N\}.
\end{align*}
Moreover, each eigenvalue is simple and the associated normalized eigenfunctions are
\begin{align*}
    \phi_k(x) = \sqrt{2} \sin(\pi k x).
\end{align*}
\item (Periodic) We have
\begin{align*}
    \sigma(-\Delta_+) = \{ 4\pi^2 (k-1)^2 : k\in \N \}.
\end{align*}
Moreover, the lowest eigenvalue $\lambda_1 = 0$ is simple and the rest is double degenerate with the following normalized eigenfunctions
\begin{align*}
    \phi_1(x) = 1, \quad \mbox{and}\quad \phi_{k,1}(x) = \sqrt{2} \cos(2 \pi (k-1) x), \quad \phi_{k,2}(x) = \sqrt{2} \sin( 2\pi (k-1) x).
\end{align*}
\item (Anti-periodic) We have
\begin{align*}
    \sigma(-\Delta_-) = \{  \pi^2 (2k-1)^2 : k \in \N\}.
\end{align*}
Moreover, every eigenvalue is double degenerate with normalized eigenfunctions
\begin{align*}
    \phi_{k,1}(x) = \sqrt{2} \cos(\pi (2k-1) x) \quad \mbox{and}\quad \phi_{k,2}(x) = \sqrt{2} \sin(\pi (2k-1) x).
\end{align*}
\end{enumerate}
\end{proposition}

\begin{proof} These results are classical and follow from straightforward calculations. 
 \end{proof}

One can now verify that Proposition~\ref{prop:1D Laplacian} is optimal in the following sense. By ordering the eigenvalues of the Laplacian in non-decreasing order and counting multiplicity $\{\lambda_k\}_{k\geq 1}$, we see that 
\begin{align*}
    \lambda_k < \lambda_{k+1},
\end{align*}
holds if and only if $k$ is odd in the periodic case and $k$ is even in the anti-periodic case. Therefore, by considering the associated $N$-particle system of electrons, i.e., the operator $-\Delta_{\pm}$ associated to $\norm{\nabla \Psi}_{\mL^2}^2$ with form domain $\mathcal{H}_N \cap \mH^1_{\pm 1}(I_N)$, we obtain explicit examples where the ground-state of the $N$-particle system is non-degenerate if and only if condition~\eqref{eq:parity condition} holds. 

In addition, note that by perturbing the Laplacian via external potentials, we can break the symmetry and make the eigenvalues of $-\Delta_{\pm}$ simple. For instance, by adding $v= \alpha \delta_x$ for $\alpha>0$ small enough we can make the first eigenvalue of $-\Delta_-$ or the second eigenvalue of $-\Delta_+$ become simple. In particular, considering the $N$-particle version of such operators, we find systems where the ground-state is non-degenerate even though the condition in~\eqref{eq:parity condition} is not fulfilled. This justifies Remark~\ref{rem:parity condition}.

\bigskip

\end{document}